\documentclass[10pt]{amsart}
\usepackage[T1]{fontenc}
\usepackage{lmodern}
\usepackage{geometry}
\usepackage[latin1] {inputenc}
\usepackage{amsmath}
\usepackage{amsfonts, amssymb, textcomp}
\usepackage[colorlinks=flase, linkcolor=red,urlcolor=green, citecolor=blue]{hyperref}
\usepackage{subeqnarray}

\usepackage{latexsym}
\usepackage{fancyhdr}
\usepackage{longtable}
\usepackage{amsmath, amssymb}
\usepackage{graphicx}

\numberwithin{equation}{section}

\theoremstyle{plain}
\newtheorem{theorem}{Theorem}[section]
\newtheorem{lemma}[theorem]{Lemma}
\newtheorem{proposition}[theorem]{Proposition}
\newtheorem{corollary}[theorem]{Corollary}
\theoremstyle{definition}
\newtheorem{definition}[theorem]{Definition}

\newtheorem{remark}[theorem]{Remark}

\newcommand{\RR}{\mathbb{R}}
\newcommand{\CC}{\mathbb{C}}
\newcommand{\NN}{\mathbb{N}}

\newcommand{\supp}{\operatorname{supp}}
\newcommand{\id}{\operatorname{id}}

\let\on=\operatorname

\newenvironment{demo}[1]{\par\smallskip\noindent{\bf #1.}}{\par\smallskip}

\makeatletter
\@namedef{subjclassname@2020}{%
  \textup{2020} Mathematics Subject Classification}
\makeatother

\title[On convolved weight matrices]
{On convolved weight matrices and local solvability with controlled loss of regularity}

\author[G.~Schindl]{Gerhard Schindl}

\address{G.~Schindl: Fakult\"at f\"ur Mathematik, Universit\"at Wien, Oskar-Morgenstern-Platz~1, A-1090 Wien, Austria.}
\email{gerhard.schindl@univie.ac.at}

\begin{document}

\begin{abstract}
We introduce the convolution between abstractly given anisotropic weight matrices and investigate properties of this new construction. Further, we apply the knowledge to the particular but interesting case when the (isotropic) matrices are associated with weight functions in the sense of Braun-Meise-Taylor and show how the convolution of associated weight matrices modifies the underlying weight functions. Finally, we give a concrete application of the convolution when studying the concept of local solvability for a certain hyperbolic PDE. Indeed, the convolution allows to treat in a natural way a mixed setting; i.e. having a controlled loss of regularity expressed in terms of two, in general different, weight sequences.
\end{abstract}

\thanks{This research was funded in whole or in part by the Austrian Science Fund (FWF) 10.55776/PAT9445424}
\date{\today}

\maketitle

\section{Introduction}
Weighted spaces appear in different contexts and fields in Mathematics and play important roles in applications in various directions. For example, we mention \emph{ultradifferentiable and ultraholomorphic classes,} see \cite{Komatsu73}, \cite{BraunMeiseTaylor90}, \cite{BonetMeiseMelikhov07}, \cite{compositionpaper}, and \cite{Thilliezdivision}, \cite{optimalflat23}; generalized \emph{Gelfand-Shilov classes,} see \cite{nuclearglobal2}, \cite{anisolog}; and \emph{weighted spaces of entire functions,} see \cite{Bonet2022survey}, \cite{weightedentireinclusion1}. In each instance we also refer to the list of citations in these works.\vspace{6pt}

Classically, one can find two approaches: either using a \emph{weight sequence} $\mathbf{M}=(M_p)_{p\in\NN}$ or a \emph{weight function} $\omega:[0,+\infty)\rightarrow[0,+\infty)$. In general, the weight sequence and weight function setting are mutually distinct, see \cite{BonetMeiseMelikhov07}, \cite[Sect. 5]{compositionpaper}, and \cite[Sect. 5]{GelfandShilovincl}. More recently, even \emph{weight matrices} or \emph{weight function systems} have been considered and these methods allow to generalize the original settings; see e.g. \cite{compositionpaper}, \cite{equalitymixedOregular}, \cite{anisolog}, \cite{ultradifferentiablecomparison}. More precisely, to each weight function $\omega$ (in the sense of Braun-Meise-Taylor) one can associate a weight matrix $\mathcal{M}_{\omega}:=\{\mathbf{W}^{(\ell)}: \ell>0\}$; see Section \ref{weightfctsection} for precise definitions and citations. In all these settings growth and regularity assumptions on the weights are required and unavoidable but, however, similar growth conditions appear and hence the frameworks are related from this point of view.\vspace{6pt}

It is a natural idea to generate new weights from given ones, for concrete applications it is frequently unavoidable to deal with such constructions and modifications. For weight sequences constructions are quite explicit and then it is interesting to understand the corresponding actions in the weight function and weight matrix setting, especially how these operations are interacting between $\omega$ and $\mathcal{M}_{\omega}$. For example, one can consider the \emph{point-wise product} $\mathbf{M}\cdot\mathbf{N}$ resp. \emph{quotient sequence} $\frac{\mathbf{M}}{\mathbf{N}}$ when $\mathbf{M}$, $\mathbf{N}$ are given and this has been studied in detail in the recent works \cite{genLegendreconj}, \cite{genLegendreconjBMT}.

The main aim of this article is to investigate the \emph{convolution} $\star$ between weight sequences and weight matrices, to show how the convolution of associated weight matrices is effecting the underlying weight functions and to give a concrete application in the matrix setting. The convolution between weight sequences has already been introduced in  \cite[$(3.15)$]{Komatsu73} and it was used and studied in \cite[Sect. 4]{weightedentireinclusion1} in the weighted entire setting and in \cite[Sect. 5]{optimalflat23} in the ultraholomorphic framework. In the latter work it turned out that even the convolution between standard weight sequences can become involved. More precisely, our goals are:

\begin{itemize}
\item[$(I)$] Define $\mathcal{M}\star\mathcal{N}$ between (abstractly given) weight matrices $\mathcal{M}:=\{\mathbf{M}^{(\iota)}: \iota>0\}$ and $\mathcal{N}:=\{\mathbf{N}^{(\iota)}: \iota>0\}$ (see Section \ref{matrixsection}), and then study the meaning resp. the counterpart of $\star$ in the weight function setting: One can ask if the matrix $\mathcal{M}_{\tau}\star\mathcal{M}_{\sigma}$ is associated with a weight function $\omega$ and how $\omega$, $\sigma$, and $\tau$ are related.

\item[$(II)$] Give the basic definitions and prove consequences and growth relations between weight matrices in the general anisotropic setting; i.e. $\mathbf{M}^{(\iota)}\in\RR_{>0}^{\NN^d}$ with $d=1,2,\dots$ arbitrary. Note that usually in higher dimensions one has $M_{|\alpha|}$ for $\alpha\in\NN^d$ and functions $\omega$ are extended radially to whole $\RR^d$.

\item[$(III)$] Study an explicit application of the convolution in a new direction; more precisely we deal with the notion of \emph{local solvability} in the theory of linear partial differential operators. The goal is to extend and generalize a result in \cite{OliaroPopivanov06} concerning a certain hyperbolic PDE from Gevrey sequences to more general weight sequences and to treat a \emph{mixed setting} with a controlled loss of regularity which is a natural notion within the weight matrix framework.
\end{itemize}

Due to the systematic and abstract flavour of $(I)$ and $(II)$ one can expect further applications of the shown results and techniques when dealing with weighted spaces (even in the anisotropic case). $(III)$ illustrates in a new direction the importance of the convolution since it allows to extend naturally the notion of local solvability to a mixed setting; i.e. when having a controlled loss of regularity in the sense that one deals with two (in general different) sequences $\mathbf{M}$ and $\mathbf{N}$ (cf. Definition \ref{locsovdef}) and this notion is becoming crucial when working with (associated) weight matrices.\vspace{6pt}

The paper is structured as follows: In Section \ref{weightsequencesection} we gather all relevant notation and conditions for weight sequences, functions, and matrices. Section \ref{mainconvolvsection} is dedicated to the detailed and systematic study of the convolution between sequences and matrices and in the main results Theorem \ref{convolvedmatrixthm} and Corollaries \ref{convolvedmatrixthmcor}, \ref{convolvedmatrixthmcor1} we investigate the role of the convolution of (associated) weight matrices and see that $\star$ corresponds to the \emph{point-wise sum $+$} between weight functions (in the sense of Braun-Meise-Taylor). Finally, in Section \ref{hyperbolicsection} we treat $(III)$ and refer to the main statement Theorem \ref{Thm12igeneral} and its variants Theorems \ref{Thm12igeneral1}, \ref{Thm12igeneral3}, and \ref{Thm12igeneral2}. First, we prove an abstract result for certain abstractly given weight sequences and then focus on associated weight matrices and hence investigate the Braun-Meise-Taylor setting in detail (see also the example constructed in Proposition \ref{counterexample}).\vspace{6pt}

\textbf{Acknowledgements.} The author thanks Chiara Boiti (Univ. di Ferrara) and Alessandro Oliaro (Univ. di Torino) for helpful discussions during the preparation of this article, in particular concerning the results contained in Section \ref{hyperbolicsection}, and Stefan F\"{u}rd\"{o}s (Univ. of Vienna) for drawing the author's attention to the paper \cite{OliaroPopivanov06}.\vspace{6pt}

\textbf{Declarations.} The author confirms that there are no competing interests and no data have been generated during writing this article.

\section{Weight sequences and weight matrices}\label{weightsequencesection}
\subsection{Basic notation}
Write $\NN=\{0,1,2,\dots\}$, $\NN_{>0}:=\NN\backslash\{0\}$ and, occasionally, we set $\RR_{>0}:=(0,+\infty)$. $\mathcal{E}$ denotes the class of smooth functions. For given multi-indices $\alpha=(\alpha_1,\dots,\alpha_d),\beta=(\beta_1,\dots,\beta_d)\in\NN^d$ we write $\alpha\le\beta$ if $\alpha_i\le\beta_i$ for all $1\le i\le d$ and set as usual $|\alpha|:=\sum_{j=1}^d\alpha_j$. On the other hand, for $t\in\RR^d$ the expression $|t|$ refers to the usual Euclidean norm on $\RR^d$.

Let $t=(t_1,\dots,t_d)\in\RR^d$, $\alpha=(\alpha_1,\dots,\alpha_d)\in\NN^d$ and $h>0$, then $t^{\alpha}:=\prod_{i=1}^dt_i^{\alpha_i}$ and $(ht)^{\alpha}=(ht_1,\dots,ht_d)^{\alpha}=\prod_{i=1}^d(ht_i)^{\alpha_i}=h^{\alpha_1+\dots+\alpha_d}t^{\alpha}=h^{|\alpha|}t^{\alpha}$.

\subsection{Anisotropic weight sequences}\label{anisosection}
Let $\mathbf{M}=(M_{\alpha})_{\alpha\in\NN^d}\in\RR_{>0}^{\NN^d}$ be given and hence such a sequence can be viewed as a map $\mathbf{M}:\NN^d\rightarrow(0,+\infty)^d$, $\alpha\mapsto M_{\alpha}$. We call $\mathbf{M}$ \emph{isotropic} if $M_{\alpha}=M_{\beta}$ for all $\alpha,\beta\in\NN^d$ satisfying $|\alpha|=|\beta|$, i.e. if $\mathbf{M}$ is only depending on the length (norm) of the multi-indices, and otherwise we say that $\mathbf{M}$ is \emph{anisotropic.} Hence in the isotropic case the crucial (growth) information is purely encoded in terms of the sequence $\mathbf{M}^r:=(M^r_p)_{p\in\NN}$ given by $M^r_p:=M_{\alpha}$, $|\alpha|=p$, and we simply write $\mathbf{M}=(M_{|\alpha|})_{\alpha\in\NN^d}$. Indeed, when $d=1$ then each sequence is isotropic.

Let $\mathbf{M}=(M_{\alpha})_{\alpha\in\NN^d}$, then analogously as in the isotropic situation in \cite[$(2.3)$]{regularnew} we put
\begin{equation}\label{liminfcond}
\mathbf{M}_{\iota}:=\liminf_{|\alpha|\rightarrow+\infty}\left(\frac{M_{\alpha}}{M_0}\right)^{1/|\alpha|}=\liminf_{|\alpha|\rightarrow+\infty}(M_{\alpha})^{1/|\alpha|}.
\end{equation}

For any $c>0$ set
\begin{equation}\label{Mcequ}
\mathbf{M}^c:=(c^{|\alpha|}M_{\alpha})_{\alpha\in\NN^d};
\end{equation}
cf. the notation in \cite[Rem. 2.10 $(a)$]{weightedentireinclusion1} in the isotropic setting. For our considerations we consider the following basic set of sequences:

\begin{definition}\label{weightsequdef}
$\mathbf{M}=(M_{\alpha})_{\alpha\in\NN^d}$ is called a \emph{weight sequence} if $$(\mathbf{M}_{\iota}=)\lim_{|\alpha|\rightarrow+\infty}(M_{\alpha})^{1/|\alpha|}=+\infty.$$
\end{definition}

Let us introduce growth relations between $\mathbf{M}=(M_{\alpha})_{\alpha\in\NN^d}$, $\mathbf{N}=(N_{\alpha})_{\alpha\in\NN^d}$:

\begin{itemize}
\item[$(*)$] Write $\mathbf{M}\le\mathbf{N}$ if $M_{\alpha}\le N_{\alpha}$ for all $\alpha\in\NN^d$.

\item[$(*)$] Write $\mathbf{M}\hypertarget{preceq}{\preceq}\mathbf{N}$ if
$$\sup_{\alpha\in\NN^d\backslash\{0\}}\left(\frac{M_{\alpha}}{N_{\alpha}}\right)^{1/|\alpha|}<+\infty;$$
i.e. if $M_{\alpha}\le C^{|\alpha|+1}N_{\alpha}$ for some $C\ge 1$ and all $\alpha\in\NN^d$.

\item[$(*)$] $\mathbf{M},\mathbf{N}\in\RR_{>0}^{\NN}$ are called \emph{equivalent,} denoted by $\mathbf{M}\hypertarget{approx}{\approx}\mathbf{N}$, if $\mathbf{M}\hyperlink{preceq}{\preceq}\mathbf{N}$ and $\mathbf{N}\hyperlink{preceq}{\preceq}\mathbf{M}$ which means that
\begin{equation}\label{equivalencrelation}
0<\inf_{\alpha\in\NN^d\backslash\{0\}}\left(\frac{M_{\alpha}}{N_{\alpha}}\right)^{1/|\alpha|}\le\sup_{\alpha\in\NN^d\backslash\{0\}}\left(\frac{M_{\alpha}}{N_{\alpha}}\right)^{1/|\alpha|}<+\infty.
\end{equation}
By definition $\mathbf{M}^c\hyperlink{approx}{\approx}\mathbf{M}$ holds for any $c>0$ and any $\mathbf{M}=(M_{\alpha})_{\alpha\in\NN^d}$.

\item[$(*)$] Write $\mathbf{M}\hypertarget{mtriangle}{\vartriangleleft}\mathbf{N}$ if $\lim_{|\alpha|\rightarrow+\infty}\left(\frac{M_{\alpha}}{N_{\alpha}}\right)^{1/|\alpha|}=0$; i.e. if
\begin{equation}\label{triangleestim}
\forall\;h>0\;\exists\;C_h\ge 1\;\forall\;\alpha\in\NN^d:\;\;\;M_{\alpha}\le C_hh^{|\alpha|}N_{\alpha}.
\end{equation}
Obviously, $\mathbf{M}\hyperlink{mtriangle}{\vartriangleleft}\mathbf{N}$ implies $\mathbf{M}\hyperlink{preceq}{\preceq}\mathbf{N}$, but \hyperlink{mtriangle}{$\vartriangleleft$} is not reflexive.
\end{itemize}

Let $\mathbf{M}=(M_{\alpha})_{\alpha\in\NN^d}$ and $\mathbf{N}=(N_{\alpha})_{\alpha\in\NN^d}$ be given, then write $(\mathbf{M},\mathbf{N})_{\on{mg}}$ if $\mathbf{M}$ and $\mathbf{N}$ satisfy
\begin{equation}\label{mg}
\exists\;C\ge 1\;\forall\;\alpha,\beta\in\NN^d:\;\;\;M_{\alpha+\beta}\le C^{|\alpha|+|\beta|+1}N_{\alpha}N_{\beta}.
\end{equation}
Similarly, write $(\mathbf{M},\mathbf{N})_{\on{dc}}$ if $\mathbf{M}$ and $\mathbf{N}$ satisfy
\begin{equation}\label{mixeddc}
\exists\;C\ge 1\;\forall\;\alpha\in\NN^d\;\forall\;1\le j\le d:\;\;\;M_{\alpha+e_j}\le C^{|\alpha|+1}N_{\alpha},
\end{equation}
with $e_j$ denoting the $j$-th unit vector in $\RR^d$. Obviously, $(\mathbf{M},\mathbf{N})_{\on{mg}}$ implies $(\mathbf{M},\mathbf{N})_{\on{dc}}$; if $(\mathbf{M},\mathbf{M})_{\on{mg}}$ holds then we say that $\mathbf{M}$ satisfies \emph{moderate growth,} abbreviated by \hypertarget{mg}{$(\on{mg})$}, whereas $\mathbf{M}$ satisfies \emph{derivation closedness,} abbreviated by \hypertarget{dc}{$(\on{dc})$}, if $(\mathbf{M},\mathbf{M})_{\on{dc}}$. In case $\mathbf{M}$ is isotropic, then \hyperlink{mg}{$(\on{mg})$} is precisely the classical condition $(M.2)$ from \cite{Komatsu73} and \hyperlink{dc}{$(\on{dc})$} corresponds to $(M.2)'$ and which motivates this terminology. When $M_0=1=N_0$, then the factor $C^{|\alpha|+|\beta|+1}$ in \eqref{mg} can be replaced by $C^{|\alpha|+|\beta|}$ and concerning \eqref{mg}, \eqref{mixeddc} we refer to \cite[$(3.6)-(3.9)$]{nuclearglobal2}.\vspace{6pt}

Next, according to \cite[Def. 5.1]{anisolog} a sequence $\mathbf{M}=(M_{\alpha})_{\alpha\in\NN^d}$ is called \emph{log-convex} if there exists a convex function $F: [0,+\infty)^d\rightarrow\RR$ such that
$$\forall\;\alpha\in\NN^d:\;\;\;F(\alpha)=\log(M_{\alpha}).$$
This notion generalizes log-convexity for isotropic sequences, i.e. $(M.1)$ from \cite{Komatsu73}, which reads as follows:
$$\forall\;p\in\NN_{>0}:\;M_p^2\le M_{p-1} M_{p+1},$$
and which is equivalent to the fact that the sequence $\mu:=(\mu_p)_{p\in\NN_{>0}}$ given by $\mu_p:=\frac{M_p}{M_{p-1}}$ is non-decreasing (we also put $\mu_0:=1$). The analogous notation is used for all appearing sequences and thus the representation $\frac{M_p}{M_0}=\prod_{i=0}^p\mu_i$ holds for all $p\in\NN$. (One can also start the product with $i=1$ and use the convention that for $p=0$ the empty product gives the value $1$.) The log-convex minorant $\mathbf{M}^{\on{lc}}$ is the largest sequence (w.r.t. the order relation $\le$) among all log-convex sequences $\mathbf{N}$ satisfying $\mathbf{N}\le\mathbf{M}$; see \cite[Sect. 5]{anisolog}.\vspace{6pt}

For isotropic sequences $\mathbf{M}$ we write that $\mathbf{M}$ is \emph{normalized} if $1=M_0\le M_1\Leftrightarrow 1\le\mu_1$ and it is convenient to consider the following set of sequences
$$\hypertarget{LCset}{\mathcal{LC}}:=\{\mathbf{M}\in\RR_{>0}^{\NN}:\;\mathbf{M}\;\text{is normalized, log-convex},\;\lim_{p\rightarrow+\infty}(M_p)^{1/p}=+\infty\}.$$
Finally, an isotropic sequence is called \emph{non-quasianalytic} if
\begin{equation}\label{nqa}
\sum_{p\ge 1}\frac{1}{\mu_p}<+\infty;
\end{equation}
this is property $(M.3)'$ in \cite{Komatsu73}. For log-convex sequences it is known that \eqref{nqa} is equivalent to $\sum_{p\ge 1}\frac{1}{(M_p)^{1/p}}<+\infty$: This holds by \emph{Carleman's inequality;} for a proof see e.g. \cite[Prop. 4.1.7 \& 4.1.8]{diploma} and the citations there.

\subsection{Weight matrices}\label{matrixsection}
Let $\mathcal{I}=\RR_{>0}$ be the index set equipped with the natural order. An \emph{(anisotropic) weight matrix} $\mathcal{M}$ associated with $\mathcal{I}$ is the set $\mathcal{M}=\{\mathbf{M}^{(\iota)}=(M^{(\iota)}_{\alpha})_{\alpha\in\NN^d}: \iota\in\mathcal{I}\}$ such that $\mathbf{M}^{(\iota)}\le\mathbf{M}^{(\iota_1)}$ for all $0<\iota\le\iota_1$; see also \cite[Sect. 3]{nuclearglobal2} and \cite[Sect. 6]{anisolog}. $\mathcal{M}$ is called
\begin{itemize}
\item[$(*)$] \emph{isotropic} if each $\mathbf{M}^{(\iota)}$ is isotropic, i.e. $\mathbf{M}^{(\iota)}=(M^{(\iota)}_{|\alpha|})_{\alpha\in\NN^d}$ for each $\iota\in\mathcal{I}$;

\item[$(*)$]\emph{constant} if $\mathbf{M}^{(\iota)}\hyperlink{approx}{\approx}\mathbf{M}^{(\iota_1)}$ for all $\iota,\iota_1\in\mathcal{I}$;

\item[$(*)$] \emph{log-convex} if $\mathbf{M}^{(\iota)}$ is log-convex for all $\iota\in\mathcal{I}$;

\item[$(*)$] \emph{standard log-convex} if $\mathcal{M}$ is isotropic and $\mathbf{M}^{(\iota)}\in\hyperlink{LCset}{\mathcal{LC}}$ for all $\iota\in\mathcal{I}$;

\item[$(*)$] \emph{non-quasianalytic} if $\mathcal{M}$ is isotropic and $\mathbf{M}^{(\iota)}$ is non-quasianalytic for all $\iota\in\mathcal{I}$.
\end{itemize}

In particular, for a standard log-convex weight matrix each $\mathbf{M}^{(\iota)}$ is an isotropic weight sequence according to Definition \ref{weightsequdef}. Let us consider some growth conditions on $\mathcal{M}$ and relations between matrices $\mathcal{M}$ and $\mathcal{N}$; see \cite[Sect. 4.1 \& 4.2]{compositionpaper}:

\hypertarget{R-mg}{$(\mathcal{M}_{\{\on{mg}\}})$} \hskip1cm $\forall\;\iota\in\mathcal{I}\;\exists\;C>0\;\exists\;\iota_1\in\mathcal{I}\;\forall\;\alpha,\beta\in\NN^d:\;\;\;M^{(\iota)}_{\alpha+\beta}\le C^{|\alpha|+|\beta|+1} M^{(\iota_1)}_{\alpha}M^{(\iota_1)}_{\beta}$,\par\vskip.3cm
\hypertarget{B-mg}{$(\mathcal{M}_{(\on{mg})})$} \hskip1cm $\forall\;\iota\in\mathcal{I}\;\exists\;C>0\;\exists\;\iota_1\in\mathcal{I}\;\forall\;\alpha,\beta\in\NN^d:\;\;\;M^{(\iota_1)}_{\alpha+\beta}\le C^{|\alpha|+|\beta|+1} M^{(\iota)}_{\alpha}M^{(\iota)}_{\beta}$,\par\vskip.3cm
\hypertarget{R-L}{$(\mathcal{M}_{\{\on{L}\}})$} \hskip1cm $\forall\;\iota\in\mathcal{I}\;\forall\;h>0\;\exists\;\iota_1\in\mathcal{I}\;\exists\;C>0\;\forall\;\alpha\in\NN^d:\;\;\;h^{|\alpha|}M^{(\iota)}_{\alpha}\le CM^{(\iota_1)}_{\alpha}$,\par\vskip.3cm
\hypertarget{B-L}{$(\mathcal{M}_{(\on{L})})$} \hskip1cm $\forall\;\iota\in\mathcal{I}\;\forall\;h>0\;\exists\;\iota_1\in\mathcal{I}\;\exists\;C>0\;\forall\;\alpha\in\NN^d:\;\;\;h^{|\alpha|}M^{(\iota_1)}_{\alpha}\le CM^{(\iota)}_{\alpha}$.\par\vskip.3cm

The brackets $\{\cdot\}$ emphasize the fact that these conditions are crucial for \emph{Roumieu-type spaces} whereas $(\cdot)$ for the \emph{Beurling-type;} see e.g. \cite[Sect. 4.2]{compositionpaper}. We write $[\cdot]$ as a joint notation meaning either $\{\cdot\}$ or $(\cdot)$ and similarly for the corresponding weighted spaces; cf. \cite[p. 99]{compositionpaper}. Consequently, \hyperlink{R-mg}{$(\mathcal{M}_{\{\on{mg}\}})$} resp. \hyperlink{B-mg}{$(\mathcal{M}_{(\on{mg})})$} precisely means that for all $\iota$ there exists some $\iota_1$ such that $(\mathbf{M}^{(\iota)},\mathbf{M}^{(\iota_1)})_{\on{mg}}$ resp. $(\mathbf{M}^{(\iota_1)},\mathbf{M}^{(\iota)})_{\on{mg}}$ holds. The relevant growth relations between weight matrices are defined as follows:
$$\mathcal{M}\{\preceq\}\mathcal{N}:\;\;\;\forall\;\iota\in\mathcal{I}\;\exists\;\iota_1\in\mathcal{I}:\;\;\;\mathbf{M}^{(\iota)}\hyperlink{preceq}{\preceq}\mathbf{N}^{(\iota_1)},$$
$$\mathcal{M}(\preceq)\mathcal{N}:\;\;\;\forall\;\iota\in\mathcal{I}\;\exists\;\iota_1\in\mathcal{I}:\;\;\;\mathbf{M}^{(\iota_1)}\hyperlink{preceq}{\preceq}\mathbf{N}^{(\iota)},$$
$$\mathcal{M}\vartriangleleft\mathcal{N}:\;\;\;\forall\;\iota,\iota_1\in\mathcal{I}:\;\;\;\mathbf{M}^{(\iota)}\hyperlink{mtriangle}{\vartriangleleft}\mathbf{N}^{(\iota_1)}.$$
We call $\mathcal{M}$ and $\mathcal{N}$ \emph{R-equivalent,} written $\mathcal{M}\{\approx\}\mathcal{N}$, if $\mathcal{M}\{\preceq\}\mathcal{N}$ and $\mathcal{N}\{\preceq\}\mathcal{M}$ and \emph{B-equivalent,} denoted by $\mathcal{M}(\approx)\mathcal{N}$, if $\mathcal{M}(\preceq)\mathcal{N}$ and $\mathcal{N}(\preceq)\mathcal{M}$. Relation $\vartriangleleft$ is not reflexive. In the isotropic setting, under mild standard assumptions on the matrices (e.g. being standard log-convex) we get that these relations are characterizing the (continuous) inclusions of the weighted spaces; see \cite[Prop. 4.6]{compositionpaper}: $\mathcal{M}\{\preceq\}\mathcal{N}$ if and only if $\mathcal{E}_{\{\mathcal{M}\}}\subseteq\mathcal{E}_{\{\mathcal{N}\}}$, $\mathcal{M}(\preceq)\mathcal{N}$ if and only if $\mathcal{E}_{(\mathcal{M})}\subseteq\mathcal{E}_{(\mathcal{N})}$, $\mathcal{M}\vartriangleleft\mathcal{N}$ if and only if $\mathcal{E}_{\{\mathcal{M}\}}\subseteq\mathcal{E}_{(\mathcal{N})}$ and the sufficiency of these condition is clear even for arbitrary (anisotropic) weight matrices.

\subsection{Associated weight function}\label{assfctsection}
Let $\mathbf{M}=(M_{\alpha})_{\alpha\in\NN^d}$ be given with $\lim_{|\alpha|\rightarrow+\infty}(M_{\alpha})^{1/|\alpha|}=+\infty$ (cf. Definition \ref{weightsequdef}). Then, in the general anisotropic setting, the \emph{associated (weight) function} has been studied in detail in \cite{anisolog}, see also \cite[Rem. 1]{nuclearglobal2}. We set
\begin{equation}\label{assoweightfct}
\omega_{\mathbf{M}}(t):=\sup_{\alpha\in\NN^d}\log\frac{M_0|t^{\alpha}|}{M_{\alpha}},\;\;\;t\in\RR^d,
\end{equation}
with the conventions $0^0:=1$ and $\log(0):=-\infty$; see \cite[Sect. 5, Rem. 5.5]{anisolog}. Alternatively, for any fixed $t\in\RR^d$ we can restrict in \eqref{assoweightfct} to all multi-indices belonging to the set $\NN_t^d:=\{\alpha=(\alpha_1,\dots,\alpha_d)\in\NN^d: \alpha_j=0\;\text{if}\;t_j=0,\;j=1,\dots,d\}$ (and again using the convention $0^0:=1$). The general assumption $\lim_{|\alpha|\rightarrow+\infty}(M_{\alpha})^{1/|\alpha|}=+\infty$ ensures the fact that $\omega_{\mathbf{M}}(t)<+\infty$ for all $t\in\RR^d$. Indeed, as seen in \cite[Rem. 1]{nuclearglobal2} this growth property even characterizes this situation and so for any weight sequence according to Definition \ref{weightsequdef} the corresponding function $\omega_{\mathbf{M}}$ is well defined (i.e. $\omega_{\mathbf{M}}(t)<+\infty$ for all $t$). Note that $\omega_{\mathbf{M}}(t)\le\omega_{\mathbf{M}}(s)$ if $t,s\in(0,+\infty)^d$ with $t\le s$ meaning that $t_j\le s_j$ for all $j=1,\dots,d$ (see \cite[Sect. 5]{anisolog}). Moreover, it is clear that $\lim_{|t|\rightarrow+\infty}\omega_{\mathbf{M}}(t)=+\infty$ since as $|t|\rightarrow+\infty$ at least one coordinate $t_j$ has to tend to infinity.

A crucial achievement in \cite{anisolog} has been to establish a connection between $\omega_{\mathbf{M}}$ and the \emph{log-convex minorant} $\mathbf{M}^{\on{lc}}$ of $\mathbf{M}$ analogously to the isotropic (resp. one-dimensional) case treated in \cite[Chapitre I]{mandelbrojtbook} or \cite{regularnew}. In \cite[Thm. 5.2 \& Rem. 5.5]{anisolog} it has been shown that
\begin{equation}\label{assoweightfctinv}
M^{\on{lc}}_{\alpha}=M_0\sup_{s\in(0,+\infty)^d}\frac{s^{\alpha}}{\exp(\omega_{\mathbf{M}}(s))},\;\;\;\forall\;\alpha\in\NN^d,
\end{equation}
which generalizes the known formula in the isotropic case from \cite[Chapitre I]{Mandelbrojt40}.\vspace{6pt}

By definition the following relation between $\mathbf{M}^c$ (recall \eqref{Mcequ}) and $\mathbf{M}$ is immediate; see \cite[Rem. 2.10, $(2.17)$]{weightedentireinclusion1} in the isotropic setting:
\begin{equation}\label{217entire}
	\forall\;c>0\;\forall\;t\in\RR^d:\;\;\;\omega_{\mathbf{M}^c}(t)=\omega_{\mathbf{M}}\left(\frac{t}{c}\right),
\end{equation}
with the notation $\frac{t}{c}=(t_1/c,\dots,t_d/c)$ when $t=(t_1,\dots,t_d)$. Note that:
\begin{itemize}
	\item[$(*)$] $\mathbf{M}^c_0=\mathbf{M}_0$ for all $c>0$;
	
	\item[$(*)$] $\lim_{|\alpha|\rightarrow+\infty}(M_{\alpha})^{1/|\alpha|}=+\infty$ if and only if $\lim_{|\alpha|\rightarrow+\infty}(M^c_{\alpha})^{1/|\alpha|}=+\infty$ for some/any $c>0$;
	
	\item[$(*)$] $\mathbf{M}$ is log-convex if and only if some/any $\mathbf{M}^c$ is so.
\end{itemize}

Let $\mathbf{M}$ be an isotropic and log-convex weight sequence, then introduce the \emph{counting function}
\begin{equation}\label{counting}
\Sigma_{\mathbf{M}}(t):=|\{p\ge 1:\;\;\;\mu_p\le t\}|.
\end{equation}
By the known \emph{integral representation formula,} see \cite[1.8. III]{mandelbrojtbook}, \cite[$(3.11)$]{Komatsu73} and also \cite[Lemma 2.5]{regularnew}, one has
\begin{equation}\label{assointrepr}
\forall\;t\ge 0:\;\;\;\omega_{\mathbf{M}}(t)=\int_0^t\frac{\Sigma_{\mathbf{M}}(u)}{u}du=\int_{\mu_1}^t\frac{\Sigma_{\mathbf{M}}(u)}{u}du.
\end{equation}

Now we turn back to the general anisotropic setting. By definition it is immediate that $\mathbf{M}\le\mathbf{N}$ implies $\omega_{\mathbf{N}}(t)\le\omega_{\mathbf{M}}(t)$ for all $t\in\RR^d$ and we characterize the more general growth relations between (anisotropic) weight sequences in terms of their corresponding associated weight functions.

\begin{lemma}\label{inclusioncharacterizationlemma}
Let $\mathbf{M},\mathbf{N}\in\RR_{>0}^{\NN^d}$ be given such that $\lim_{|\alpha|\rightarrow+\infty}(M_{\alpha})^{1/|\alpha|}=+\infty$ (i.e. $\mathbf{M}$ is a weight sequence according to Definition \ref{weightsequdef}).
\begin{itemize}
\item[$(i)$] $\mathbf{M}\hyperlink{preceq}{\preceq}\mathbf{N}$ implies
\begin{equation}\label{inclusioncharacterizationlemmaequ}
\exists\;C,h\ge 1\;\forall\;t\in\RR^d:\;\;\;\omega_{\mathbf{N}}(t)\le\omega_{\mathbf{M}}(ht)+C,
\end{equation}
and $\mathbf{M}\hyperlink{mtriangle}{\vartriangleleft}\mathbf{N}$ implies
\begin{equation}\label{inclusioncharacterizationlemmaequ1}
\forall\;h>0\;\exists\;C_h\ge 1\;\forall\;t\in\RR^d:\;\;\;\omega_{\mathbf{N}}(t)\le\omega_{\mathbf{M}}(ht)+C_h.
\end{equation}
\item[$(ii)$] If $\mathbf{M}$ is also log-convex, then \eqref{inclusioncharacterizationlemmaequ} implies $\mathbf{M}\hyperlink{preceq}{\preceq}\mathbf{N}$ whereas \eqref{inclusioncharacterizationlemmaequ1} implies $\mathbf{M}\hyperlink{mtriangle}{\vartriangleleft}\mathbf{N}$.
\end{itemize}
Moreover, the assertions in $(i)$ and $(ii)$ imply $\lim_{|\alpha|\rightarrow+\infty}(N_{\alpha})^{1/|\alpha|}=+\infty$ too (equivalently $\omega_{\mathbf{N}}$ is well defined).
\end{lemma}

\emph{Note:} \eqref{inclusioncharacterizationlemmaequ1} is precisely relation $\omega_{\mathbf{M}}\vartriangleleft_{\mathfrak{c}}\omega_{\mathbf{N}}$ in \cite[Sect. 2.3]{ultradifferentiablecomparison}.

\demo{Proof}
$(i)$ First, by definition relation $\mathbf{M}\hyperlink{preceq}{\preceq}\mathbf{N}$, and hence $\mathbf{M}\hyperlink{mtriangle}{\vartriangleleft}\mathbf{N}$, implies $\lim_{|\alpha|\rightarrow+\infty}(N_{\alpha})^{1/|\alpha|}=+\infty$. Concerning \eqref{inclusioncharacterizationlemmaequ}, by assumption $\mathbf{M}\hyperlink{preceq}{\preceq}\mathbf{N}$ one gets $M_{\alpha}\le Ch^{|\alpha|}N_{\alpha}$ for some $C,h\ge 1$ and all $\alpha\in\NN^d$ and so $\frac{N_0|t^{\alpha}|}{N_{\alpha}}\le\frac{C N_0}{M_0}\frac{M_0|(ht)^{\alpha}|}{M_{\alpha}}$ for all $t\in\RR^d$, $\alpha\in\NN^d$. Applying $\log$ yields the conclusion and \eqref{inclusioncharacterizationlemmaequ1} follows analogously since in this situation even for each $h>0$ there exists $C_h>0$ such that the above estimate is valid.\vspace{6pt}

$(ii)$ We combine \eqref{inclusioncharacterizationlemmaequ} and \eqref{assoweightfctinv} and get for all $\alpha\in\NN^d$:
\begin{align*}
N_{\alpha}&\ge N^{\on{lc}}_{\alpha}=N_0\sup_{s\in(0,+\infty)^d}\frac{s^{\alpha}}{\exp(\omega_{\mathbf{N}}(s))}\ge N_0e^{-C}\sup_{s\in(0,+\infty)^d}\frac{s^{\alpha}}{\exp(\omega_{\mathbf{M}}(hs))}
\\&
=N_0e^{-C}\sup_{u\in(0,+\infty)^d}\frac{(h^{-1}u)^{\alpha}}{\exp(\omega_{\mathbf{M}}(u))}=\frac{N_0}{M_0e^C}h^{-|\alpha|}M_0\sup_{u\in(0,+\infty)^d}\frac{u^{\alpha}}{\exp(\omega_{\mathbf{M}}(u))}
\\&
=\frac{N_0}{M_0e^C}h^{-|\alpha|}M_{\alpha},
\end{align*}
which gives the conclusion. We have applied the substitution $u=hs$, i.e. set $u_j:=hs_j$ for all $j=1,\dots,d$ and the second part follows analogously. Finally, note that both \eqref{inclusioncharacterizationlemmaequ} and \eqref{inclusioncharacterizationlemmaequ1} yield, in particular, that $\omega_{\mathbf{N}}(t)<+\infty$ for all $t\in\RR^d$ since $\omega_{\mathbf{M}}(t)<+\infty$ for all $t\in\RR^d$ by assumption $\lim_{|\alpha|\rightarrow+\infty}(M_{\alpha})^{1/|\alpha|}=+\infty$. Hence $\lim_{|\alpha|\rightarrow+\infty}(N_{\alpha})^{1/|\alpha|}=+\infty$ follows by \cite[Rem. 1]{nuclearglobal2}.
\qed\enddemo

\subsection{Weight functions and associated weight matrices}\label{weightfctsection}
We introduce weights in the sense of \emph{Braun-Meise-Taylor (BMT-weights for short).} Let $\omega:[0,+\infty)\rightarrow[0,+\infty)$ be continuous, non-decreasing, $\omega(0)=0$ and $\lim_{t\rightarrow+\infty}\omega(t)=+\infty$. If $\omega$ satisfies in addition $\omega(t)=0$ for all $t\in[0,1]$, then $\omega$ is called \emph{normalized} and this can always be assumed w.l.o.g. For convenience we write that $\omega$ has $\hypertarget{om0}{(\omega_0)}$ if it satisfies all these assumptions; see e.g. \cite[Sect. 2.1]{index} and \cite[Sect. 2.2]{sectorialextensions}.\vspace{6pt}

We consider the following (standard) conditions; these abbreviations have already been used in ~\cite{dissertation}.

\begin{itemize}
\item[\hypertarget{om1}{$(\omega_1)$}] $\omega(2t)=O(\omega(t))$ as $t\rightarrow+\infty$; i.e. $\exists\;L\ge 1\;\forall\;t\ge 0:\;\;\;\omega(2t)\le L\omega(t)+L$.

\item[\hypertarget{om3}{$(\omega_3)$}] $\log(t)=o(\omega(t))$ as $t\rightarrow+\infty$.
	
\item[\hypertarget{om4}{$(\omega_4)$}] $\varphi_{\omega}:t\mapsto\omega(e^t)$ is a convex function on $\RR$.

\item[\hypertarget{om6}{$(\omega_6)$}] $\exists\;H\ge 1\;\forall\;t\ge 0:\;\;\;2\omega(t)\le\omega(Ht)+H$.

\item[\hypertarget{omnq}{$(\omega_{\text{nq}})$}] $\int_1^{\infty}\frac{\omega(t)}{t^2}dt<+\infty$.
\end{itemize}

For any $s>0$ we write $\id^{1/s}: t\mapsto t^{1/s}$ \emph{(Gevrey weights).} In the literature one can find different assumptions for BMT-weights (cf. \cite[Def. 1.4]{ultradifferentiablecomparison}), however, some of the conditions are basic and for convenience we introduce the set
$$\hypertarget{omset0}{\mathcal{W}_0}:=\{\omega:[0,\infty)\rightarrow[0,\infty): \omega\;\text{has}\;\hyperlink{om0}{(\omega_0)},\hyperlink{om3}{(\omega_3)},\hyperlink{om4}{(\omega_4)}\}.$$
In the forthcoming \hyperlink{omset0}{$\mathcal{W}_0$} is understood to be the set of (normalized) weight functions in the sense of Braun-Meise-Taylor. Then, for any $\omega\in\hyperlink{omset0}{\mathcal{W}_0}$ we define the \emph{Legendre-Fenchel-Young-conjugate} of $\varphi_{\omega}$ by
\begin{equation}\label{legendreconjugate}
	\varphi^{*}_{\omega}(x):=\sup\{x y-\varphi_{\omega}(y): y\ge 0\},\;\;\;x\ge 0,
\end{equation}
with the following properties, see e.g. \cite[Rem. 1.3, Lemma 1.5]{BraunMeiseTaylor90} and \cite[Sect. 2.9]{ultradifferentiablecomparison}: It is convex and non-decreasing, $\varphi^{*}_{\omega}(0)=0$, $\varphi^{**}_{\omega}=\varphi_{\omega}$, $\lim_{x\rightarrow+\infty}\frac{x}{\varphi^{*}_{\omega}(x)}=0$ and finally $x\mapsto\frac{\varphi_{\omega}(x)}{x}$ and $x\mapsto\frac{\varphi^{*}_{\omega}(x)}{x}$ are non-decreasing on $(0,+\infty)$. Note that by normalization we can extend the supremum in \eqref{legendreconjugate} from $y\ge 0$ to $y\in\RR$ without changing the value of $\varphi^{*}_{\omega}(x)$ for any given $x\ge 0$.\vspace{6pt}

If $\omega\in\hyperlink{omset0}{\mathcal{W}_0}$ satisfies \hyperlink{omnq}{$(\omega_{\on{nq}})$}, then $\omega$ is called \emph{non-quasianalytic} and \emph{quasianalytic} else.\vspace{6pt}

We introduce now crucial growth relations; see also \cite[Sect. 2.3]{ultradifferentiablecomparison}. Let $\sigma,\tau: [0,+\infty)\rightarrow[0,+\infty)$ be functions satisfying $\lim_{t\rightarrow+\infty}\sigma(t)=+\infty=\lim_{t\rightarrow+\infty}\tau(t)$.

\begin{itemize}
\item[$(*)$] Write $\sigma\hypertarget{ompreceq}{\preceq}\tau$ if
\begin{equation}\label{bigOrelation}
\tau(t)=O(\sigma(t))\;\text{as}\;t\rightarrow+\infty.	
\end{equation}
\item[$(*)$] $\sigma$ and $\tau$ are called \emph{equivalent,} denoted by $\sigma\hypertarget{sim}{\sim}\tau$, if $\sigma\hyperlink{ompreceq}{\preceq}\tau$ and $\tau\hyperlink{ompreceq}{\preceq}\sigma$.

\item[$(*)$] Write $\sigma\hypertarget{omvartriangle}{\vartriangleleft}\tau$ if
\begin{equation}\label{smallOrelation}
	\tau(t)=o(\sigma(t))\;\text{as}\;t\rightarrow+\infty,	
\end{equation}
and note that \hyperlink{omvartriangle}{$\vartriangleleft$} is not reflexive.
\end{itemize}
The analogous relations can be considered for functions $\sigma,\tau:[0,+\infty)^d\rightarrow[0,+\infty)$ as $|t|\rightarrow+\infty$ when $\lim_{|t|\rightarrow+\infty}\sigma(t)=+\infty=\lim_{|t|\rightarrow+\infty}\tau(t)$.

We mention the following known result, see e.g. \cite[Lemma 2.8]{testfunctioncharacterization} resp. \cite[Lemma 2.4]{sectorialextensions} and the references mentioned in the proofs there.

\begin{lemma}\label{assoweightomega0}
Let $\mathbf{M}\in\hyperlink{LCset}{\mathcal{LC}}$. Then $\omega_{\mathbf{M}}\in\hyperlink{omset0}{\mathcal{W}_0}$. Furthermore, \hyperlink{om6}{$(\omega_6)$} holds for $\omega_{\mathbf{M}}$ if and only if $\mathbf{M}$ satisfies \hyperlink{mg}{$(\on{mg})$}.
\end{lemma}

Now let us summarize some facts which are shown in \cite[Sect. 5]{compositionpaper} and are needed in this work; all properties listed below are valid for $\omega\in\hyperlink{omset0}{\mathcal{W}_0}$ except \eqref{newexpabsorb} for which \hyperlink{om1}{$(\omega_1)$} is crucial. More basic properties and conditions for abstractly given (isotropic) weight matrices $\mathcal{M}$ can be found e.g. in \cite[Sect. 4]{compositionpaper}.

\begin{itemize}
	\item[$(i)$] The idea was that to each $\omega\in\hyperlink{omset0}{\mathcal{W}_0}$ one can associate an isotropic weight matrix $\mathcal{M}_{\omega}:=\{\mathbf{W}^{(\ell)}=(W^{(\ell)}_p)_{p\in\NN}: \ell>0\}$ by\vspace{6pt}
	
	\centerline{$W^{(\ell)}_p:=\exp\left(\frac{1}{\ell}\varphi^{*}_{\omega}(\ell p)\right)$.}\vspace{6pt}
	
	We have that $\mathbf{W}^{(\ell)}\in\hyperlink{LCset}{\mathcal{LC}}$ for each $\ell>0$, i.e. $\mathcal{M}_{\omega}$ is \emph{standard log-convex,} and we even have the stronger order relation
	\begin{equation*}\label{quotientorderequ}
		\forall\;\ell_2\ge\ell_1>0:\;\;\;\vartheta^{(\ell_1)}\le\vartheta^{(\ell_2)},
	\end{equation*}
	with $\vartheta^{(\ell)}$ denoting the corresponding sequence of quotients; see \cite[Sect. 2.5]{whitneyextensionweightmatrix}.
	
	\item[$(ii)$] $\mathcal{M}_{\omega}$ satisfies
	\begin{equation}\label{newmoderategrowth}
		\forall\;\ell>0\;\forall\;p,q\in\NN:\;\;\;W^{(\ell)}_{p+q}\le W^{(2\ell)}_pW^{(2\ell)}_q;
	\end{equation}
	i.e. $(\mathbf{W}^{(\ell)},\mathbf{W}^{(2\ell)})_{\on{mg}}$ for all $\ell>0$ and thus both \hyperlink{R-mg}{$(\mathcal{M}_{\{\on{mg}\}})$} and \hyperlink{B-mg}{$(\mathcal{M}_{(\on{mg})})$} hold.
	
	\item[$(iii)$] \hyperlink{om6}{$(\omega_6)$} for $\omega$ holds if and only if some/each $\mathbf{W}^{(\ell)}$ satisfies \hyperlink{mg}{$(\on{mg})$} if and only if $\mathbf{W}^{(\ell_1)}\hyperlink{approx}{\approx}\mathbf{W}^{(\ell_2)}$ for each $\ell_1,\ell_2>0$. Thus \hyperlink{om6}{$(\omega_6)$} is characterizing the case when $\mathcal{M}_{\omega}$ is \emph{constant.}
	
	\item[$(iv)$] In case $\omega$ has in addition \hyperlink{om1}{$(\omega_1)$}, then $\mathcal{M}_{\omega}$ also satisfies
	\begin{equation}\label{newexpabsorb}
		\forall\;h\ge 1\;\exists\;d\ge 1\;\forall\;\ell>0\;\exists\;D\ge 1\;\forall\;p\in\NN:\;\;\;h^pW^{(\ell)}_p\le D W^{(d\ell)}_p.
	\end{equation}
 This estimate implies \hyperlink{R-L}{$(\mathcal{M}_{\{\on{L}\}})$} and \hyperlink{B-L}{$(\mathcal{M}_{(\on{L})})$} and is crucial for verifying the equality $\mathcal{E}_{[\mathcal{M}_{\omega}]}=\mathcal{E}_{[\omega]}$ (as l.c.v.s.). Concerning the definition of these spaces we refer to \cite{BraunMeiseTaylor90} and \cite{compositionpaper}; recall that the notation $[\cdot]$ is the convention meaning either the \emph{Roumieu case} $\{\cdot\}$ or the \emph{Beurling case} $(\cdot)$. The same equalities hold for analogously defined weighted spaces since the implication showing \eqref{newexpabsorb} purely involves properties of the defining weight function; see \cite[Lemma 5.9, Thm. 5.14]{compositionpaper}.
	
	\item[$(v)$] We have $\omega\hyperlink{sim}{\sim}\omega_{\mathbf{W}^{(\ell)}}$ for each $\ell>0$, more precisely
	\begin{equation}\label{goodequivalenceclassic}
		\forall\;\ell>0\,\,\exists\,D_{\ell}>0\;\forall\;t\ge 0:\;\;\;\ell\omega_{\mathbf{W}^{(\ell)}}(t)\le\omega(t)\le 2\ell\omega_{\mathbf{W}^{(\ell)}}(t)+D_{\ell};
	\end{equation}
	for a proof see \cite[Theorem 4.0.3, Lemma 5.1.3]{dissertation}, \cite[Lemma 5.7]{compositionpaper} and also \cite[Lemma 2.5]{sectorialextensions}. Note that, on the one hand, for proving \eqref{goodequivalenceclassic} the convexity condition \hyperlink{om4}{$(\omega_4)$} is indispensable but, on the other hand, \hyperlink{om4}{$(\omega_4)$} is only required for the second estimate in \eqref{goodequivalenceclassic}.

\item[$(vi)$] It holds that $\mathcal{M}_{\omega}$ is non-quasianalytic if and only if $\omega$ satisfies \hyperlink{omnq}{$(\omega_{\on{nq}})$}: For this use the fact that this condition is preserved under equivalence, \eqref{goodequivalenceclassic} and \cite[Lemma 4.1]{Komatsu73} applied to (each) $\mathbf{W}^{(\ell)}$; see also \cite[Cor. 4.8]{testfunctioncharacterization}.
\end{itemize}
The growth relations are crucial and characterizing the inclusions in the ultradifferentiable setting; more precisely by \cite[Lemma 5.16 \& Cor. 5.17]{compositionpaper} for $\sigma,\tau\in\hyperlink{omset0}{\mathcal{W}_0}$ satisfying \hyperlink{om1}{$(\omega_1)$} one has $\mathcal{E}_{[\sigma]}\subseteq\mathcal{E}_{[\tau]}$ if and only if $\sigma\hyperlink{ompreceq}{\preceq}\tau$ and $\mathcal{E}_{\{\sigma\}}\subseteq\mathcal{E}_{(\tau)}$ if and only if $\sigma\hyperlink{omvartriangle}{\vartriangleleft}\tau$. Note that for the sufficiency it is enough to assume $\sigma,\tau\in\hyperlink{omset0}{\mathcal{W}_0}$ and this part also holds for analogously defined weighted spaces by the similar seminorms; see the proof of \cite[Lemma 5.16]{compositionpaper}.\vspace{6pt}

In view of Lemma \ref{assoweightomega0}, for any $\mathbf{M}\in\hyperlink{LCset}{\mathcal{LC}}$ we can naturally consider $\mathcal{M}_{\omega_{\mathbf{M}}}:=\{\mathbf{M}^{(\iota)}: \iota>0\}$ and in this case $\mathbf{M}=\mathbf{M}^{(1)}$; see e.g. \cite[$(2.13)$, p. 407]{subaddlike}.

\section{Convolved weights}\label{mainconvolvsection}
We introduce convolved (anisotropic) sequences and weight matrices, the corresponding counter-part in the Braun-Meise-Taylor weight function setting and, finally, study the effects of the convolution on the associated weighted (ultradifferentiable) function classes.

\subsection{Convolved sequences}\label{convolvesection}
Let $\mathbf{M},\mathbf{N}\in\RR_{>0}^{\NN^d}$, then define the \emph{convolved sequence} $\mathbf{M}\star\mathbf{N}=(M\star N_\beta)_{\beta\in\NN^d}$ by
\begin{equation}\label{convolvesequ}
M\star N_{\beta}:=\min_{0\le\alpha\le\beta}M_{\alpha}N_{\beta-\alpha},\;\,\;\beta\in\NN^d.
\end{equation}
For the corresponding definition in the isotropic case see \cite[$(3.15)$]{Komatsu73} and \cite[Sect. 4.2]{weightedentireinclusion1}. Note that for simplicity we avoid the notation $(M\star N)_{\beta}$ and write $M\star N_{\beta}$ instead. Next let us summarize some immediate consequences; see also \cite[Sect. 4.2]{weightedentireinclusion1}:

\begin{itemize}
\item[$(a)$] $\mathbf{M}\star\mathbf{N}=\mathbf{N}\star\mathbf{M}$ is valid; i.e. $\star$ is commutative.

\item[$(b)$] $M\star N_0=M_0N_0$ and thus $M_0=1=N_0$ implies $M\star N_0=1$.

\item[$(c)$] More generally, $M\star N_{\alpha}\le\min\{M_0N_{\alpha},M_{\alpha}N_0\}$ for all $\alpha\in\NN^d$ and so $$\mathbf{M}\star\mathbf{N}\le\max\{M_0,N_0\}\min\{\mathbf{M},\mathbf{N}\}.$$

Consequently, if $M_0=1=N_0$ then $\mathbf{M}\star\mathbf{N}\le\min\{\mathbf{M},\mathbf{N}\}$ and if $\mathbf{M}\le\mathbf{N}$, then $\mathbf{M}\star\mathbf{N}\le N_0\mathbf{M}$.

\item[$(d)$] By definition $(\mathbf{M},\mathbf{N})_{\on{mg}}$ amounts to have $\mathbf{M}\hyperlink{preceq}{\preceq}\mathbf{N}\star\mathbf{N}$ and $\mathbf{M}$ satisfies \hyperlink{mg}{$(\on{mg})$} if and only if $\mathbf{M}\hyperlink{preceq}{\preceq}\mathbf{M}\star\mathbf{M}$, i.e. if and only if $\mathbf{M}\hyperlink{approx}{\approx}\mathbf{M}\star\mathbf{M}$. Even more generally one has $\mathbf{L}\hyperlink{preceq}{\preceq}\mathbf{M}\star\mathbf{N}$ if and only if
    \begin{equation*}\label{mixedmg}
    \exists\;C\ge 1\;\forall\;\alpha,\beta\in\NN^d:\;\;\;L_{\alpha+\beta}\le C^{|\alpha|+|\beta|+1}M_{\alpha}N_{\beta},
    \end{equation*}
    and in this condition the ``relative growth'' between $\mathbf{M}$ and $\mathbf{N}$ can behave irregular, i.e. neither $\mathbf{M}\hyperlink{preceq}{\preceq}\mathbf{N}$ nor $\mathbf{N}\hyperlink{preceq}{\preceq}\mathbf{M}$ holds.
\end{itemize}

For the next properties we focus on the isotropic case:

\begin{itemize}
\item[$(e)$] Let $\mathbf{M}$, $\mathbf{N}$ be log-convex, then by \cite[Lemma 3.5]{Komatsu73} the corresponding quotient sequence $\mu\star\nu$ is obtained when rearranging resp. ordering the elements of the sequences $\mu$ and $\nu$ in the order of growth: For this fact recall that by $(b)$
    $$M\star N_q=\min_{0\le p\le q}M_0\mu_1\cdots\mu_p\cdot N_0\nu_1\cdots\nu_{q-p}=M\star N_0\min_{0\le p\le q}\mu_1\cdots\mu_p\cdot\nu_1\cdots\nu_{q-p},$$
    with the convention that the empty product corresponds to the value $1$.

    This identity implies the fact that the convolution of log-convex (weight) sequences is again a log-convex (weight) sequence according to Definition \ref{weightsequdef} and $\mathbf{M}\star\mathbf{N}\in\hyperlink{LCset}{\mathcal{LC}}$ provided that $\mathbf{M},\mathbf{N}\in\hyperlink{LCset}{\mathcal{LC}}$.

\item[$(f)$] Let $\mathbf{M},\mathbf{N}$ be log-convex weight sequences, then by $(e)$ and the definition of the counting function in \eqref{counting} we get
\begin{equation}\label{countingadditive}
\forall\;t\ge 0:\;\;\;\Sigma_{\mathbf{M}\star\mathbf{N}}(t)=\Sigma_{\mathbf{M}}(t)+\Sigma_{\mathbf{N}}(t),
\end{equation}
and by \eqref{assointrepr} this identity implies
\begin{equation}\label{assofctadditive}
\forall\;t\ge 0:\;\;\;\omega_{\mathbf{M}\star\mathbf{N}}(t)=\omega_{\mathbf{M}}(t)+\omega_{\mathbf{N}}(t).
\end{equation}
\end{itemize}

\begin{remark}\label{commentfrem}
Note that for $(f)$ it is natural to assume $\mathbf{M}_{\iota}=+\infty=\mathbf{N}_{\iota}$ which is ensured when dealing with weight sequences in the sense of Definition \ref{weightsequdef}, recall also \eqref{liminfcond}. But this does not follow automatically by the assumptions on the sequences considered in \cite[Lemma 3.5]{Komatsu73}: Log-convexity only implies $\mathbf{M}_{\iota},\mathbf{N}_{\iota}>0$; see \cite[Sect. 2.2, Lemma 2.1]{regularnew} and also \cite[Rem. 4.1]{conjugateweightfunction}. However, note that $\omega_{\mathbf{M}}(t)=+\infty$ for all $t>\mathbf{M}_{\iota}$ if $\mathbf{M}_{\iota}<+\infty$, see \cite[Lemma 2.2]{regularnew} and also \cite[Lemmas 2.4, 2.5 \& 2.6]{regularnew}, and then the equalities in \eqref{countingadditive} and \eqref{assofctadditive} are (formally) still valid for all $t\ge 0$ when setting both sides to be $+\infty$ for all $t>\min\{\mathbf{M}_{\iota},\mathbf{N}_{\iota}\}=\mathbf{M}\star\mathbf{N}_{\iota}$ and this equality follows immediately by \cite[Sect. 2.2, Lemma 2.1]{regularnew} and $(e)$ above.
\end{remark}

The analogue of \eqref{countingadditive} in the anisotropic setting is unclear since the appropriate generalization of the counting function $\Sigma_{\mathbf{M}}$ is required. On the other hand, we can directly prove that \eqref{assofctadditive} is valid in the anisotropic setting as well and the idea applies also to the isotropic setting. Therefore, \eqref{assofctadditive} even holds without having assumption log-convexity for $\mathbf{M}$ and/or $\mathbf{N}$.

\begin{lemma}\label{anisoassofctadditivelemma}
Let $\mathbf{M},\mathbf{N}\in\RR_{>0}^{\NN^d}$ be given weight sequences (see Definition \ref{weightsequdef}). Then
\begin{equation}\label{anisoassofctadditivelemmaequ}
\forall\;t\in\RR^d:\;\;\;\omega_{\mathbf{M}\star\mathbf{N}}(t)=\omega_{\mathbf{M}}(t)+\omega_{\mathbf{N}}(t).
\end{equation}
\end{lemma}

\demo{Proof}
In order to show \eqref{anisoassofctadditivelemmaequ} we apply the exponential function to it and in view of \eqref{assoweightfct} this precisely means
$$\forall\;t\in\RR^d:\;\;\;\sup_{\alpha\in\NN^d}\frac{M_0|t^{\alpha}|}{M_{\alpha}}\sup_{\beta\in\NN^d}\frac{N_0|t^{\beta}|}{N_{\beta}}=\sup_{\gamma\in\NN^d}\frac{M\star N_0|t^{\gamma}|}{M\star N_{\gamma}}.$$
So fix $t\in\RR^d$. First, for any $\alpha,\beta\in\NN^d$ set $\gamma_{\alpha,\beta}:=\alpha+\beta$ and get $t^{\gamma_{\alpha,\beta}}=t^{\alpha}t^{\beta}$ and $M\star N_{\gamma_{\alpha,\beta}}=\min_{0\le\delta\le\gamma_{\alpha,\beta}}M_{\delta}N_{\gamma_{\alpha,\beta}-\delta}\le M_{\alpha}N_{\beta}$ when choosing $\delta:=\alpha$. Recall also that by $(b)$ above $M\star N_0=M_0N_0$ is valid. Altogether, this implies $$\exp(\omega_{\mathbf{M}\star\mathbf{N}}(t))=\sup_{\gamma\in\NN^d}\frac{M\star N_0|t^{\gamma}|}{M\star N_{\gamma}}\ge\sup_{\alpha,\beta\in\NN^d}\frac{M_0|t^{\alpha}|}{M_{\alpha}}\frac{N_0|t^{\beta}|}{N_{\beta}}=\exp(\omega_{\mathbf{M}}(t))\exp(\omega_{\mathbf{N}}(t)).$$
Conversely, let $\gamma\in\NN^d$ be given and then we can find $\alpha_{\gamma},\beta_{\gamma}\in\NN^d$ with $\alpha_{\gamma},\beta_{\gamma}\le\gamma$ and such that $M\star N_{\gamma}=\min_{0\le\delta\le\gamma}M_{\delta}N_{\gamma-\delta}=M_{\alpha_{\gamma}}N_{\beta_{\gamma}}$. Therefore, $\alpha_{\gamma}+\beta_{\gamma}=\gamma$ and because $t^{\gamma}=t^{\alpha_{\gamma}}t^{\beta_{\gamma}}$ and $M\star N_0=M_0N_0$ we get
$$\sup_{\gamma\in\NN^d}\frac{M\star N_0|t^\gamma|}{M\star N_\gamma}\le\sup_{\alpha_{\gamma}\in\NN^d}\frac{M_0|t^{\alpha_{\gamma}}|}{M_{\alpha_{\gamma}}}\sup_{\beta_{\gamma}\in\NN^d}\frac{N_0|t^{\beta_{\gamma}}|}{N_{\beta_{\gamma}}}\le\sup_{\alpha\in\NN^d}\frac{M_0|t^\alpha|}{M_{\alpha}}\sup_{\beta\in\NN^d}\frac{N_0|t^\beta|}{N_{\beta}}.$$
\qed\enddemo

\begin{corollary}\label{anisoassofctadditivecor}
Let $\mathbf{L},\mathbf{N}\in\RR_{>0}^{\NN^d}$ be given weight sequences (see Definition \ref{weightsequdef}). Then we get:
\begin{itemize}
\item[$(i)$] $(\mathbf{L},\mathbf{N})_{\on{mg}}$ implies
\begin{equation}\label{om6aniso}
\exists\;H\ge 1\;\forall\;t\in\RR^d:\;\;\;2\omega_{\mathbf{N}}(t)\le\omega_{\mathbf{L}}(Ht)+H.
\end{equation}
\item[$(ii)$] If $\mathbf{L}$ is in addition log-convex, then \eqref{om6aniso} implies $(\mathbf{L},\mathbf{N})_{\on{mg}}$.
\end{itemize}
\end{corollary}

\demo{Proof}
$(i)$ $(\mathbf{L},\mathbf{N})_{\on{mg}}$ means $\mathbf{L}\hyperlink{preceq}{\preceq}\mathbf{N}\star\mathbf{N}$, see $(d)$ above, and then by $(i)$ in Lemma \ref{inclusioncharacterizationlemma} (see \eqref{inclusioncharacterizationlemmaequ}) we get
$$\exists\;C,h\ge 1\;\forall\;t\in\RR^d:\;\;\;\omega_{\mathbf{N}\star\mathbf{N}}(t)\le\omega_{\mathbf{L}}(ht)+C.$$
Now apply \eqref{anisoassofctadditivelemmaequ} to $\mathbf{M}=\mathbf{N}$ and this gives \eqref{om6aniso} with $H:=\max\{h,C\}$.\vspace{6pt}

$(ii)$ Conversely, we combine \eqref{om6aniso} with \eqref{assoweightfctinv} and \eqref{anisoassofctadditivelemmaequ} to obtain for all $\alpha\in\NN^d$:
\begin{align*}
L_{\alpha}&=L_0\sup_{s\in(0,+\infty)^d}\frac{s^{\alpha}}{\exp(\omega_{\mathbf{L}}(s))}\le e^HL_0\sup_{s\in(0,+\infty)^d}\frac{s^{\alpha}}{\exp(2\omega_{\mathbf{N}}(H^{-1}s))}
\\&
=\frac{e^HL_0}{N\star N_0}N\star N_0\sup_{u\in(0,+\infty)^d}\frac{(uH)^{\alpha}}{\exp(\omega_{\mathbf{N}\star\mathbf{N}}(u))}=\frac{e^HL_0H^{|\alpha|}}{N\star N_0}(N\star N)^{\on{lc}}_{\alpha}\le\frac{e^HL_0H^{|\alpha|}}{N\star N_0}N\star N_{\alpha}.
\end{align*}
This estimate means that $\mathbf{L}\hyperlink{preceq}{\preceq}\mathbf{N}\star\mathbf{N}$, i.e. $(\mathbf{L},\mathbf{N})_{\on{mg}}$.
\qed\enddemo

\begin{proposition}\label{basicprop}
We get the following properties for the convolution $\star$:
\begin{itemize}
\item[$(a)$] Let $\mathbf{M},\mathbf{N},\mathbf{R},\mathbf{S}\in\RR_{>0}^{\NN^d}$. If $\mathbf{M}\le\mathbf{R}$ and $\mathbf{N}\le\mathbf{S}$, then $\mathbf{M}\star\mathbf{N}\le\mathbf{R}\star\mathbf{S}$. If $\mathbf{M}\hyperlink{preceq}{\preceq}\mathbf{R}$ and $\mathbf{N}\hyperlink{preceq}{\preceq}\mathbf{S}$, then $\mathbf{M}\star\mathbf{N}\hyperlink{preceq}{\preceq}\mathbf{R}\star\mathbf{S}$ and, consequently, $\mathbf{M}\hyperlink{approx}{\approx}\mathbf{R}$ and $\mathbf{N}\hyperlink{approx}{\approx}\mathbf{S}$ implies $\mathbf{M}\star\mathbf{N}\hyperlink{approx}{\approx}\mathbf{R}\star\mathbf{S}$.

    And if $\mathbf{M}\hyperlink{mtriangle}{\vartriangleleft}\mathbf{R}$ and $\mathbf{N}\hyperlink{mtriangle}{\vartriangleleft}\mathbf{S}$, then $\mathbf{M}\star\mathbf{N}\hyperlink{mtriangle}{\vartriangleleft}\mathbf{R}\star\mathbf{S}$.

\item[$(b)$] Let $\mathbf{M},\mathbf{N},\mathbf{R},\mathbf{S}$ be weight sequences. If
\begin{equation}\label{mixedom1}
\omega_{\mathbf{M}}(2t)=O(\omega_{\mathbf{R}}(t)),\hspace{15pt}\omega_{\mathbf{N}}(2t)=O(\omega_{\mathbf{S}}(t)),\;\;\;|t|\rightarrow+\infty,
\end{equation}
then $\omega_{\mathbf{M}\star\mathbf{N}}(2t)=O(\omega_{\mathbf{R}\star\mathbf{S}}(t))$ as $|t|\rightarrow+\infty$.

\item[$(c)$] Let $\mathbf{M},\mathbf{N},\mathbf{R},\mathbf{S}$ be weight sequences. If
\begin{equation}\label{mixedom6}
\exists\;H\ge 1\;\forall\;t\in\RR^d:\;\;\;2\omega_{\mathbf{M}}(t)\le\omega_{\mathbf{R}}(Ht)+H,\hspace{15pt}2\omega_{\mathbf{N}}(t)\le\omega_{\mathbf{S}}(Ht)+H,
\end{equation}
then
$$\forall\;t\in\RR^d:\;\;\;2\omega_{\mathbf{M}\star\mathbf{N}}(t)\le\omega_{\mathbf{R}\star\mathbf{S}}(2Ht)+2H.$$

\item[$(d)$] Let $\mathbf{M},\mathbf{N},\mathbf{R},\mathbf{S}$ be weight sequences. When $\omega_{\mathbf{M}}\hyperlink{ompreceq}{\preceq}\omega_{\mathbf{R}}$ and $\omega_{\mathbf{N}}\hyperlink{ompreceq}{\preceq}\omega_{\mathbf{S}}$, then $\omega_{\mathbf{M}\star\mathbf{N}}\hyperlink{ompreceq}{\preceq}\omega_{\mathbf{R}\star\mathbf{S}}$. The analogous implication holds for the relation \hyperlink{omvartriangle}{$\vartriangleleft$}.

    Consequently, if $\omega_{\mathbf{M}}\hyperlink{sim}{\sim}\omega_{\mathbf{R}}$ and $\omega_{\mathbf{N}}\hyperlink{sim}{\sim}\omega_{\mathbf{S}}$ then $\omega_{\mathbf{M}\star\mathbf{N}}\hyperlink{sim}{\sim}\omega_{\mathbf{R}\star\mathbf{S}}$.
\end{itemize}
\end{proposition}

\emph{Note:} $(a)$ and $(d)$ show that $\star$ preserves growth relations between sequences and associated weight functions; $(b)$ and $(c)$ imply that $\star$ also preserves mixed \hyperlink{om1}{$(\omega_1)$}- and mixed \hyperlink{om6}{$(\omega_6)$}-type conditions in the general anisotropic case (recall Section \ref{weightfctsection}).

\demo{Proof}
$(a)$ The first part is immediate by definition. Concerning the second one, by assumption there exist $C_i,h_i>0$, $i=1,2$, such that $M_{\alpha}\le C_1h_1^{|\alpha|}R_{\alpha}$ and $N_{\alpha}\le C_2h_2^{|\alpha|}S_{\alpha}$ for all $\alpha\in\NN^d$. Then set $C:=\max\{C_1,C_2\}$, $h:=\max\{h_1,h_2\}$ and get for all $\beta\in\NN^d$ that $M\star N_{\beta}=\min_{0\le\alpha\le\beta}M_{\alpha}N_{\beta-\alpha}\le C^2h^{|\beta|}\min_{0\le\alpha\le\beta}R_{\alpha}S_{\beta-\alpha}=C^2h^{|\beta|}R\star S_{\beta}$.

Finally, in the third part note that in this case even for all $h>0$ (small) there exists $C_h>0$ (large) such that $M_{\alpha}\le C_hh^{|\alpha|}R_{\alpha}$ and $N_{\alpha}\le C_hh^{|\alpha|}S_{\alpha}$ for all $\alpha\in\NN^d$ and one concludes as before.\vspace{6pt}

$(b)$ \eqref{mixedom1} yields $\omega_{\mathbf{M}}(2t)\le L_1\omega_{\mathbf{R}}(t)+L_1$ and $\omega_{\mathbf{N}}(2t)\le L_2\omega_{\mathbf{S}}(t)+L_2$ for some $L_i\ge 1$, $i=1,2$, and all $t\in\RR^d$ and recall the notation $2t=(2t_1,\dots,2t_d)$. Then the conclusion follows by applying Lemma \ref{anisoassofctadditivelemma} (see also \eqref{assofctadditive} in the log-convex isotropic case):
$$\omega_{\mathbf{M}\star\mathbf{N}}(2t)=\omega_{\mathbf{M}}(2t)+\omega_{\mathbf{N}}(2t)\le L_1\omega_{\mathbf{R}}(t)+L_2\omega_{\mathbf{S}}(t)+L_1+L_2\le L\omega_{\mathbf{R}\star\mathbf{S}}(t)+L;$$
indeed it suffices to choose $L:=2\max\{L_1,L_2\}$.\vspace{6pt}

$(c)$ In view of Lemma \ref{anisoassofctadditivelemma} (see also \eqref{assofctadditive} in the log-convex isotropic case) we get for all $t\in\RR^d$:
$$2\omega_{\mathbf{M}\star\mathbf{N}}(t)=2\omega_{\mathbf{M}}(t)+2\omega_{\mathbf{N}}(t)\le\omega_{\mathbf{R}}(Ht)+\omega_{\mathbf{S}}(Ht)+2H=\omega_{\mathbf{R}\star\mathbf{S}}(Ht)+2H,$$
which shows the assertion.\vspace{6pt}

$(d)$ By assumption $\omega_{\mathbf{R}}(t)\le C_1\omega_{\mathbf{M}}(t)+C_1$ and $\omega_{\mathbf{S}}(t)\le C_2\omega_{\mathbf{N}}(t)+C_2$ for some $C_i\ge 1$ and all $t\in\RR^d$. Set $C:=\max\{C_1,C_2\}$ and then Lemma \ref{anisoassofctadditivelemma} implies for all $t\in\RR^d$:
$$\omega_{\mathbf{R}\star\mathbf{S}}(t)=\omega_{\mathbf{R}}(t)+\omega_{\mathbf{S}}(t)\le C\omega_{\mathbf{M}}(t)+C\omega_{\mathbf{N}}(t)+2C=C\omega_{\mathbf{M}\star\mathbf{N}}(t)+2C.$$
The second part follows analogously since here even for all $c>0$ there exists $C\ge 1$ such that $\omega_{\mathbf{R}}(t)\le c\omega_{\mathbf{M}}(t)+C$ and $\omega_{\mathbf{S}}(t)\le c\omega_{\mathbf{N}}(t)+C$ for all $t\in\RR^d$.
\qed\enddemo

The next result shows that (in the isotropic setting) the convolution preserves non-quasianalyticity.

\begin{lemma}\label{nqaconvolvelemma}
Let $\mathbf{M},\mathbf{N}\in\RR_{>0}^{\NN}$ be log-convex weight sequences. Then $\mathbf{M}\star\mathbf{N}$ is a log-convex weight sequence and non-quasianalytic if and only if both $\mathbf{M}$ and $\mathbf{N}$ are non-quasianalytic.
\end{lemma}

\demo{Proof}
First, by $(e)$ above we see that $\mathbf{M}\star\mathbf{N}$ is also a log-convex weight sequence. Second, by \cite[Lemma 4.1]{Komatsu73} it holds that $\mathbf{M}$ resp. $\mathbf{N}$ is non-quasianalytic if and only if $\omega_{\mathbf{M}}$ resp. $\omega_{\mathbf{N}}$ is non-quasianalytic; i.e. it satisfies \hyperlink{omnq}{$(\omega_{\on{nq}})$}. Indeed, the proof of this result holds for any log-convex (isotropic) weight sequence. Then note that condition \hyperlink{omnq}{$(\omega_{\on{nq}})$} is obviously preserved under taking the pointwise sum of functions (cf. \cite[Sect. 3.9]{dissertation}).

Consequently, if both $\mathbf{M}$ and $\mathbf{N}$ are non-quasianalytic, then via Lemma \ref{anisoassofctadditivelemma} (resp. \eqref{assofctadditive}) we get that $\omega_{\mathbf{N}}+\omega_{\mathbf{N}}=\omega_{\mathbf{M}\star\mathbf{N}}$ is non-quasianalytic and hence \cite[Lemma 4.1]{Komatsu73} gives that $\mathbf{M}\star\mathbf{N}$ is non-quasianalytic, too.

On the other hand, if either $\mathbf{M}$ or $\mathbf{N}$ is quasianalytic, then by $(c)$ above we infer $(M\star N_p)^{1/p}\le\max\{(M_0)^{1/p},(N_0)^{1/p}\}\min\{(M_p)^{1/p},(N_p)^{1/p}\}$ for all $p\in\NN_{>0}$ and so, in any case, $\sum_{p\ge 1}\frac{1}{(M\star N_p)^{1/p}}=+\infty$. Thus $\mathbf{M}\star\mathbf{N}$ is quasianalytic.
\qed\enddemo

\subsection{Convolved weight matrices}\label{convolvedmatrixsection}
Let $\mathcal{M}:=\{\mathbf{M}^{(\iota)}: \iota\in\mathcal{I}\}$ and $\mathcal{N}:=\{\mathbf{N}^{(\iota)}: \iota\in\mathcal{I}\}$ be weight matrices, then inspired by \eqref{convolvesequ} we introduce the corresponding \emph{convolved weight matrix} by
\begin{equation}\label{convolvedmatrix}
\mathcal{M}\star\mathcal{N}:=\{\mathbf{M}^{(\iota)}\star\mathbf{N}^{(\iota)}: \iota\in\mathcal{I}\}.
\end{equation}
Indeed, $\mathcal{M}\star\mathcal{N}$ is a matrix according to Section \ref{matrixsection} since for any $\iota\le\iota_1$ we get $\mathbf{M}^{(\iota)}\le\mathbf{M}^{(\iota_1)}$, $\mathbf{N}^{(\iota)}\le\mathbf{N}^{(\iota_1)}$ and $(a)$ in Proposition \ref{basicprop} yields $\mathbf{M}^{(\iota)}\star\mathbf{N}^{(\iota)}\le\mathbf{M}^{(\iota_1)}\star\mathbf{N}^{(\iota_1)}$. And $(a)$ in Proposition \ref{basicprop} also implies that if both $\mathcal{M}$ and $\mathcal{N}$ are constant then $\mathcal{M}\star\mathcal{N}$, too. Moreover, clearly $\mathcal{M}\star\mathcal{N}=\mathcal{N}\star\mathcal{M}$ (cf. $(a)$ in Section \ref{convolvesection}) and if both $\mathcal{M}$ and $\mathcal{N}$ are (isotropic and) standard log-convex, then $\mathcal{M}\star\mathcal{N}$ too by $(e)$ in Section \ref{convolvesection}. Finally, if both $\mathcal{M}$ and $\mathcal{N}$ are isotropic and log-convex then $\mathcal{M}\star\mathcal{N}$ too and in this case via Lemma \ref{nqaconvolvelemma} the matrix $\mathcal{M}\star\mathcal{N}$ is non-quasianalytic if and only if both $\mathcal{M}$ and $\mathcal{N}$ are so.

As a special case consider $\mathcal{N}=\{\mathbf{N}\}$ for some (fixed) $\mathbf{N}\in\RR_{>0}^{\NN^d}$, then write
\begin{equation}\label{convolvedconstmatrix}
\mathcal{M}\star\mathbf{N}:=\{\mathbf{M}^{(\iota)}\star\mathbf{N}: \iota\in\mathcal{I}\}.
\end{equation}
Indeed, if $\mathcal{N}:=\{\mathbf{N}^{(\iota)}: \iota\in\mathcal{I}\}$ is constant then by $(a)$ in Proposition \ref{basicprop} the matrices $\mathcal{M}\star\mathcal{N}$ and $\mathcal{M}\star\mathbf{N}^{(\iota_0)}$ are $R$- and $B$-equivalent for any $\iota_0\in\mathcal{I}$. Assume that all $\mathbf{M}^{(\iota)}$, $\mathbf{N}^{(\iota)}$ are weight sequences according to Definition \ref{weightsequdef}, then we set (cf. Lemma \ref{anisoassofctadditivelemma} and \eqref{assofctadditive} in the isotropic setting):
\begin{equation}\label{convolvedfunctionmatrix}
\omega_{\mathcal{M}\star\mathcal{N}}:=\{\omega_{\mathbf{M}^{(\iota)}\star\mathbf{N}^{(\iota)}}: \iota\in\mathcal{I}\}=\{\omega_{\mathbf{M}^{(\iota)}}+\omega_{\mathbf{N}^{(\iota)}}: \iota\in\mathcal{I}\}=:\omega_{\mathcal{M}}+\omega_{\mathcal{N}}.
\end{equation}
Recall that under this assumption on the sequences also each $\omega_{\mathbf{M}^{(\iota)}\star\mathbf{N}^{(\iota)}}$ is well defined and, in particular, the sets in \eqref{convolvedfunctionmatrix} can be introduced for arbitrary standard log-convex matrices $\mathcal{M}$ and $\mathcal{N}$. Note that here one has the order $\omega_{\mathbf{M}^{(\iota)}\star\mathbf{N}^{(\iota)}}\ge\omega_{\mathbf{M}^{(\iota_1)}\star\mathbf{N}^{(\iota_1)}}$ for all $0<\iota\le\iota_1$; recall Section \ref{assfctsection}. Concerning the definition of \eqref{convolvedfunctionmatrix}, i.e. of a \emph{weight function matrix} we refer to \cite[Sect. 2.6 \& 2.7]{equalitymixedOregular} and \cite[Sect. 2.7]{testfunctioncharacterization} for the isotropic case, which means that the weight functions under consideration are \emph{radial.} More recently in \cite{ultradifferentiablecomparison} also non-radial weights have been considered (anisotropic case) but in this recent work we have reversed the order relation of the indices.

When $\sigma,\tau\in\hyperlink{omset0}{\mathcal{W}_0}$, then $\sigma+\tau\in\hyperlink{omset0}{\mathcal{W}_0}$ and when both weights satisfy in addition \hyperlink{om1}{$(\omega_1)$}, then $\sigma+\tau$ too; see \cite[Sect. 5, p. 118]{compositionpaper} and \cite[Sect. 3.9]{dissertation}: pointwise addition yields the algebraic structure of an abelian semigroup on the set \hyperlink{omset0}{$\mathcal{W}_0$} and preserves \hyperlink{om1}{$(\omega_1)$}. We prove now the main result when dealing with convolved weight matrices.

\begin{theorem}\label{convolvedmatrixthm}
Let $\mathcal{M}:=\{\mathbf{M}^{(\iota)}: \iota\in\mathcal{I}\}$ and $\mathcal{N}:=\{\mathbf{N}^{(\iota)}: \iota\in\mathcal{I}\}$ be standard log-convex, let $\mathcal{M}\star\mathcal{N}$ and $\omega_{\mathcal{M}\star\mathcal{N}}$ be the convolved matrices from \eqref{convolvedmatrix} and \eqref{convolvedfunctionmatrix}.
\begin{itemize}
\item[$(I)$] Assume the following Roumieu-type assertions for both $\mathcal{M}$ and $\mathcal{N}$:
\begin{itemize}
\item[$(a)$]
\begin{equation}\label{mixedom1roumvar}
\exists\;r>1\;\forall\;\iota\in\mathcal{I}\;\exists\;\iota_1\in\mathcal{I}\;\exists\;L\in\NN_{>0}:\;\;\;\liminf_{j\rightarrow+\infty}\frac{(W^{(\iota_1)}_{Lj})^{\frac{1}{Lj}}}{(W^{(\iota)}_j)^{\frac{1}{j}}}>r;
\end{equation}
see \cite[$(3.5)$]{equalitymixedOregular},

\item[$(b)$] and \hyperlink{R-mg}{$(\mathcal{M}_{\{\on{mg}\}})$}.
\end{itemize}
Then $\mathcal{M}\star\mathcal{N}$ satisfies \hyperlink{R-mg}{$(\mathcal{M}_{\{\on{mg}\}})$}, \eqref{mixedom1roumvar} and, as l.c.v.s.,
\begin{equation}\label{convolvedmatrixthmroumequ}
\mathcal{E}_{\{\mathcal{M}\star\mathcal{N}\}}=\mathcal{E}_{\{\omega_{\mathcal{M}\star\mathcal{N}}\}}.
\end{equation}

\item[$(II)$] Assume the following Beurling-type assertions for both matrices:
\begin{itemize}
\item[$(a)$]
\begin{equation}\label{mixedom1beurvar}
\exists\;r>1\;\forall\;\iota\in\mathcal{I}\;\exists\;\iota_1\in\mathcal{I}\;\exists\;L\in\NN_{>0}:\;\;\;\liminf_{j\rightarrow+\infty}\frac{(W^{(\iota)}_{Lj})^{\frac{1}{Lj}}}{(W^{(\iota_1)}_j)^{\frac{1}{j}}}>r;
\end{equation}
see \cite[$(3.11)$]{equalitymixedOregular},

\item[$(b)$] and \hyperlink{B-mg}{$(\mathcal{M}_{(\on{mg})})$}.
\end{itemize}
Then $\mathcal{M}\star\mathcal{N}$ satisfies \hyperlink{B-mg}{$(\mathcal{M}_{(\on{mg})})$}, \eqref{mixedom1beurvar} and, as l.c.v.s.,
\begin{equation}\label{convolvedmatrixthmbeurequ}
\mathcal{E}_{(\mathcal{M}\star\mathcal{N})}=\mathcal{E}_{(\omega_{\mathcal{M}\star\mathcal{N}})}.
\end{equation}
When in addition also $\omega_{\mathbf{M}^{(\iota)}}\hyperlink{omsim}{\sim}\omega_{\mathbf{M}^{(\iota_1)}}$ and $\omega_{\mathbf{N}^{(\iota)}}\hyperlink{omsim}{\sim}\omega_{\mathbf{N}^{(\iota_1)}}$ for all $\iota,\iota_1>0$, then
\begin{equation}\label{convolvedmatrixthmequextended}
\forall\;\iota>0:\;\;\;\mathcal{E}_{[\mathcal{M}\star\mathcal{N}]}=\mathcal{E}_{[\omega_{\mathbf{M}^{(\iota)}\star\mathbf{N}^{(\iota)}}]}=\mathcal{E}_{[\omega_{\mathbf{M}^{(\iota)}}+\omega_{\mathbf{N}^{(\iota)}}]}.
\end{equation}
Recall again that $[\cdot]$ is the joint notation meaning either $\{\cdot\}$ or $(\cdot)$.
\end{itemize}
\end{theorem}

\demo{Proof}
$(I)$ By \cite[Thm. 3.2 $(I)(ii)\Leftrightarrow(iii)$]{equalitymixedOregular} we get that \eqref{mixedom1roumvar} (for $\mathcal{M}$) is equivalent to
\begin{equation}\label{mixedom1roum}
\forall\;\iota\in\mathcal{I}\;\exists\;\iota_1\in\mathcal{I}:\;\;\;\limsup_{t\rightarrow+\infty}\frac{\omega_{\mathbf{M}^{(\iota_1)}}(2t)}{\omega_{\mathbf{M}^{(\iota)}}(t)}<+\infty;
\end{equation}
i.e. \cite[$(3.4)$]{equalitymixedOregular} and similarly for the matrix $\mathcal{N}$ with an index $\iota_2$ subject to given (and fixed) $\iota$. Hence the assumption gives \eqref{mixedom1} for $\mathbf{R}=\mathbf{M}^{(\iota)}$, $\mathbf{M}=\mathbf{M}^{(\iota_1)}$, $\mathbf{S}=\mathbf{N}^{(\iota)}$, $\mathbf{N}=\mathbf{N}^{(\iota_2)}$. Set $\overline{\iota}:=\max\{\iota_1,\iota_2\}$ and so $\mathbf{M}^{(\overline{\iota})}\ge\mathbf{M}^{(\iota_1)}$, $\mathbf{N}^{(\overline{\iota})}\ge\mathbf{N}^{(\iota_2)}$. Then, since the order relation for the associated weight functions is reversed, $(b)$ in Proposition \ref{basicprop} yields the fact that \eqref{mixedom1roum} (for both matrices) implies $\omega_{\mathbf{M}^{(\overline{\iota})}\star\mathbf{N}^{(\overline{\iota})}}(2t)=O(\omega_{\mathbf{M}^{(\iota)}\star\mathbf{N}^{(\iota)}}(t))$ as $t\rightarrow+\infty$ and hence \eqref{mixedom1roumvar} for the convolved matrix $\mathcal{M}\star\mathcal{N}$ too by applying \cite[Thm. 3.2]{equalitymixedOregular} to $\mathcal{M}\star\mathcal{N}$.

Next, by \cite[Prop. 4.1 $(i)\Leftrightarrow(ii)$]{equalitymixedOregular} (Roumieu part) it follows that \hyperlink{R-mg}{$(\mathcal{M}_{\{\on{mg}\}})$} (for $\mathcal{M}$) is equivalent to having
$$\forall\;\iota\in\mathcal{I}\;\exists\;\iota_1\in\mathcal{I}\;\exists\;H\ge 1\;\forall\;t\ge 0:\;\;\;2\omega_{\mathbf{M}^{(\iota_1)}}(t)\le\omega_{\mathbf{M}^{(\iota)}}(Ht)+H,$$
and similarly for the matrix $\mathcal{N}$ with an index $\iota_2$ (when $\iota$ is given and fixed). So \eqref{mixedom6} follows with $\mathbf{R}=\mathbf{M}^{(\iota)}$, $\mathbf{M}=\mathbf{M}^{(\iota_1)}$, $\mathbf{S}=\mathbf{N}^{(\iota)}$, $\mathbf{N}=\mathbf{N}^{(\iota_2)}$ and with $\overline{\iota}:=\max\{\iota_1,\iota_2\}$ part $(c)$ in Proposition \ref{basicprop} yields
$$\forall\;t\ge 0:\;\;\;2\omega_{\mathbf{M}^{(\overline{\iota})}\star\mathbf{N}^{(\overline{\iota})}}(t)\le\omega_{\mathbf{M}^{(\iota)}\star\mathbf{N}^{(\iota)}}(2Ht)+2H;$$
i.e. \hyperlink{R-mg}{$(\mathcal{M}_{\{\on{mg}\}})$} for $\mathcal{M}\star\mathcal{N}$ and for this recall that this convolved matrix is standard log-convex, too.

Summarizing, we are in position to apply the Roumieu part of the main result \cite[Thm. 5.4]{equalitymixedOregular} to the matrix $\mathcal{M}\star\mathcal{N}$ and this implies \eqref{convolvedmatrixthmroumequ}.\vspace{6pt}

$(II)$ The Beurling-setting follows analogously when involving \cite[Thm. 3.2 $(II)(ii)\Leftrightarrow(iii)$]{equalitymixedOregular}, the Beurling part of \cite[Prop. 4.1 $(i)\Leftrightarrow(ii)$]{equalitymixedOregular}, and finally the Beurling part of \cite[Thm. 5.4]{equalitymixedOregular}.\vspace{6pt}

Concerning the supplement \eqref{convolvedmatrixthmequextended}, the first equality holds by taking into the shown identities and $(d)$ in Proposition \ref{basicprop} which gives $\omega_{\mathbf{M}^{(\iota)}\star\mathbf{N}^{(\iota)}}\hyperlink{sim}{\sim}\omega_{\mathbf{M}^{(\iota_1)}\star\mathbf{N}^{(\iota_1)}}$ for all $\iota,\iota_1\in\mathcal{I}$. And the second equality follows directly by \eqref{assofctadditive}.
\qed\enddemo

When focusing now on associated weight matrices we obtain the following result in the BMT-weight function setting:

\begin{corollary}\label{convolvedmatrixthmcor}
Let $\sigma,\tau\in\hyperlink{omset0}{\mathcal{W}_0}$ be given with associated weight matrices $\mathcal{M}_{\sigma}:=\{\mathbf{S}^{(\ell)}: \ell>0\}$ and $\mathcal{M}_{\tau}:=\{\mathbf{T}^{(\ell)}: \ell>0\}$. Assume that both weights satisfy \hyperlink{om1}{$(\omega_1)$}, then as l.c.v.s.
\begin{equation}\label{convolvedmatrixthmcorequextended}
\forall\;\ell>0:\;\;\;\mathcal{E}_{[\mathcal{M}_{\sigma}\star\mathcal{N}_{\tau}]}=\mathcal{E}_{[\omega_{\mathbf{S}^{(\ell)}\star\mathbf{T}^{(\ell)}}]}=\mathcal{E}_{[\omega_{\mathbf{S}^{(\ell)}}+\omega_{\mathbf{T}^{(\ell)}}]}=\mathcal{E}_{[\sigma+\tau]}=\mathcal{E}_{[\mathcal{M}_{\sigma+\tau}]}.
\end{equation}
\end{corollary}

\demo{Proof}
We apply Theorem \ref{convolvedmatrixthm} to $\mathcal{M}_{\sigma}$ and $\mathcal{M}_{\tau}$: First note that by $(i)$ and $(ii)$ in Section \ref{weightfctsection} both matrices are standard log-convex and satisfy both \hyperlink{R-mg}{$(\mathcal{M}_{\{\on{mg}\}})$} and \hyperlink{B-mg}{$(\mathcal{M}_{(\on{mg})})$}(recall \eqref{newmoderategrowth}). And by \hyperlink{om1}{$(\omega_1)$} for both associated matrices \eqref{newexpabsorb} holds (recall $(iv)$ in Section \ref{weightfctsection}) and thus both \eqref{mixedom1roumvar} and \eqref{mixedom1beurvar}. More precisely, even both \hyperlink{R-L}{$(\mathcal{M}_{\{\on{L}\}})$} and \hyperlink{B-L}{$(\mathcal{M}_{(\on{L})})$} are valid by \hyperlink{om1}{$(\omega_1)$}; see also \cite[Cor. 3.3 \& Rem. 3.4]{equalitymixedOregular}.

Finally, via \eqref{goodequivalenceclassic} we have $\sigma\hyperlink{sim}{\sim}\omega_{\mathbf{S}^{(\ell)}}$ and $\tau\hyperlink{sim}{\sim}\omega_{\mathbf{T}^{(\ell)}}$ for each $\ell>0$ and so all associated weight functions (w.r.t. the particular matrix) are equivalent.

Summarizing, \eqref{convolvedmatrixthmequextended} is valid and which gives the first two equalities in \eqref{convolvedmatrixthmcorequextended}. The third ones follow again by \eqref{goodequivalenceclassic} for both matrices since this yields $\sigma+\tau\hyperlink{sim}{\sim}\omega_{\mathbf{S}^{(\ell)}}+\omega_{\mathbf{T}^{(\ell_1)}}$ for all $\ell,\ell_1>0$. Finally, the last equality holds by applying the main result \cite[Thm. 5.14 (2)]{compositionpaper} to the weight function $\sigma+\tau$ and for this recall that $\sigma+\tau\in\hyperlink{omset0}{\mathcal{W}_0}$ and satisfies \hyperlink{om1}{$(\omega_1)$}, too.
\qed\enddemo

Note that the above proof also gives that in \eqref{convolvedmatrixthmcorequextended} we can use resp. extend the chain of equalities by $\mathcal{E}_{[\omega_{\mathbf{S}^{(\ell_1)}\star\mathbf{T}^{(\ell_2)}}]}$, $\mathcal{E}_{[\omega_{\mathbf{S}^{(\ell_1)}}+\omega_{\mathbf{T}^{(\ell_2)}}]}$ with $\ell_1,\ell_2>0$ arbitrary. The next result clarifies the relation between the convolved matrix $\mathcal{M}_{\sigma}\star\mathcal{N}_{\tau}$ and $\mathcal{M}_{\sigma+\tau}$.

\begin{corollary}\label{convolvedmatrixthmcor1}
Let $\sigma,\tau\in\hyperlink{omset0}{\mathcal{W}_0}$ be given with associated weight matrices $\mathcal{M}_{\sigma}:=\{\mathbf{S}^{(\ell)}: \ell>0\}$ and $\mathcal{M}_{\tau}:=\{\mathbf{T}^{(\ell)}: \ell>0\}$.
\begin{itemize}
\item[$(i)$] $\mathcal{M}_{\sigma}\star\mathcal{N}_{\tau}$ is non-quasianalytic if and only if both $\sigma$ and $\tau$ are non-quasianalytic.

\item[$(ii)$] If both weights satisfy \hyperlink{om1}{$(\omega_1)$}, then
\begin{equation}\label{weightmatrixequivequ}
\mathcal{M}_{\sigma}\star\mathcal{N}_{\tau}\{\approx\}\mathcal{M}_{\sigma+\tau},\hspace{15pt}\mathcal{M}_{\sigma}\star\mathcal{N}_{\tau}(\approx)\mathcal{M}_{\sigma+\tau}.
\end{equation}
\end{itemize}
\end{corollary}

\emph{Note:} In view of \eqref{legendreconjugate} and the definition of the associated weight matrix a direct computation of $\mathcal{M}_{\sigma+\tau}$ in terms of $\sigma$, $\tau$ is unclear and involved.

\demo{Proof}
$(i)$ This follows by Lemma \ref{nqaconvolvelemma} and $(vi)$ in Section \ref{weightfctsection}: $\mathcal{M}_{\sigma}\star\mathcal{N}_{\tau}$ is non-quasianalytic if and only if both $\mathcal{M}_{\sigma}$ and $\mathcal{N}_{\tau}$ are non-quasianalytic which holds if and only if both $\sigma$ and $\tau$ satisfy \hyperlink{omnq}{$(\omega_{\on{nq}})$}. (Recall that taking the pointwise sum of (BMT-)weights preserves \hyperlink{omnq}{$(\omega_{\on{nq}})$}.)\vspace{6pt}

$(ii)$ \eqref{weightmatrixequivequ} is a direct consequence of the equality $\mathcal{E}_{[\mathcal{M}_{\sigma}\star\mathcal{N}_{\tau}]}=\mathcal{E}_{[\mathcal{M}_{\sigma+\tau}]}$ from \eqref{convolvedmatrixthmcorequextended} and the characterization \cite[Prop. 4.6 (1)]{compositionpaper} applied to the appearing associated weight matrices. Indeed, both $\mathcal{M}_{\sigma}\star\mathcal{N}_{\tau}$ and $\mathcal{M}_{\sigma+\tau}$ are standard log-convex and hence this result can be used.
\qed\enddemo

\begin{remark}\label{convolvedmatrixthmcorrem}
All equalities in \eqref{convolvedmatrixthmroumequ}, \eqref{convolvedmatrixthmbeurequ}, \eqref{convolvedmatrixthmequextended} and \eqref{convolvedmatrixthmcorequextended} hold for analogously defined weighted settings as well since the proofs exclusively deal with weights and their growth properties (recall also the proof of \cite[Thm. 5.4]{equalitymixedOregular}). In particular, this comment applies to ultradifferentiable test function spaces, i.e. spaces of ultradifferentiable functions having compact support and hence replacing the symbol/functor $\mathcal{E}$ by $\mathcal{D}$ (see e.g. \cite{Komatsu73}).
\end{remark}

Since these test function spaces are also crucial for the considerations in Section \ref{hyperbolicsection}, we give now the precise definition: Let $U\subseteq\RR^d$ be open and $\mathcal{M}$ be a weight matrix, then the corresponding class of test functions (resp. of functions having compact support) in $U$ is given by:
\begin{equation}\label{testfctspace}
\mathcal{D}_{[\mathcal{M}]}(U):=\{f\in\mathcal{E}_{[\mathcal{M}]}(U):\;\;\;\exists\;K\;\text{compact},\;K\subset\subset U,\;\;\;\supp(f)\subseteq K\}=\mathcal{E}_{[\mathcal{M}]}(U)\cap\mathcal{D}(U).
\end{equation}
Thus $\mathcal{D}_{\{\mathbf{G}^s\}}(U)$ (with $s>1$) is corresponding to the notation $G^s_0(U)$ in \cite{OliaroPopivanov06} and $\mathcal{D}_{[\mathbf{M}]}(U)$ corresponds to $\mathcal{D}^{[M_p]}(U)$ in \cite{Komatsu73} and in both situations one exclusively deals with isotropic weight sequences.

When now $\mathcal{M}$ is isotropic and log-convex, then $\mathcal{D}_{[\mathcal{M}]}(U)\neq\{0\}$ if and only if $\mathcal{M}$ is non-quasianalytic which follows by \cite[Sect. 4; Thm. 4.1 \& Prop. 4.7]{testfunctioncharacterization}. The characterization in terms of condition non-quasianalyticity for $\mathcal{M}$ is even established under more general assumptions on $\mathcal{M}$ since log-convexity is not required necessarily and one involves regularizations of the corresponding weight sequences in the matrix. Indeed, the proof is based on the known \emph{Denjoy-Carleman-Theorem} for the weight sequence setting; see e.g. \cite[Thm. 1.3.8]{hoermander} and \cite[Thm. 4.2]{Komatsu73}. Moreover, this general statement can be applied to $\mathcal{M}_{\omega}$ and in the BMT-weight function setting condition \hyperlink{omnq}{$(\omega_{\on{nq}})$} for $\omega$ is characterizing the non-triviality of $\mathcal{D}_{[\mathcal{M}_{\omega}]}(U)=\mathcal{D}_{[\omega]}(U)$; see \cite[Cor. 4.8]{testfunctioncharacterization} and also \cite[Sect. 2 \& 3]{BraunMeiseTaylor90} resp. \cite[Sect. 5]{compositionpaper} for the equality of the spaces.\vspace{6pt}

Thus when dealing in the above results with weighted test function spaces $\mathcal{D}_{[\mathcal{M}\star\mathcal{N}]}$ resp. $\mathcal{D}_{[\sigma+\tau]}$ then it is natural to assume that $\mathcal{M}\star\mathcal{N}$ resp. $\sigma+\tau$ is non-quasianalytic in order to ensure non-triviality of the corresponding weighted spaces. For this, in view of Lemma \ref{nqaconvolvelemma} resp. $(i)$ in Corollary \ref{convolvedmatrixthmcor1} condition non-quasianalyticity for both $\mathcal{M}$ and $\mathcal{N}$ resp. \hyperlink{omnq}{$(\omega_{\on{nq}})$} for both $\sigma$ and $\tau$ is relevant and this is consistent with the characterization in \cite[Cor. 4.8]{testfunctioncharacterization}. (Indeed, recall that by $(vi)$ in Section \ref{weightfctsection} we have that for any $\omega\in\hyperlink{omset0}{\mathcal{W}_0}$ the matrix $\mathcal{M}_{\omega}$ is non-quasianalytic if and only if $\omega$ satisfies \hyperlink{omnq}{$(\omega_{\on{nq}})$}, see also step $(b)$ in the proof of Proposition \ref{counterexample}, but in this case in general the equality $\mathcal{D}_{[\mathcal{M}_{\omega}]}(U)=\mathcal{D}_{[\omega]}(U)$ and hence the non-triviality of $\mathcal{D}_{[\omega]}(U)$ is unclear.)

\section{On mixed local solvability for an hyperbolic PDO}\label{hyperbolicsection}
We apply now the convolution in order to investigate the notion of \emph{mixed local solvability} for a certain hyperbolic PDO crucially studied in the approach in \cite{OliaroPopivanov06}; more precisely the plan is to generalize the local solvability part \cite[Thm. 1.2 $(i)$]{OliaroPopivanov06} when allowing a controlled loss of regularity. It turns out that the convolution is tailor-made to study this loss and hence this operation illustrates that the techniques in \cite{OliaroPopivanov06}, there purely exploited for the \emph{Gevrey sequences} $\mathbf{G}^s:=(p!^s)_{p\in\NN}$ with $s>1$ (and the Roumieu-type), allow for a more general statement. First, in the next section even in the anisotropic setting we prove several technical auxiliary results which can become crucial for related topics and are of independent interest.

\subsection{Technical preliminaries}\label{preliminarysection}
We introduce the following new and crucial growth relations between arbitrary $\mathbf{M},\mathbf{N},\mathbf{L}\in\RR_{>0}^{\NN^d}$; for this recall the notation from \eqref{Mcequ}:
\begin{equation}\label{PDOmainrelationgen}
\forall\;a\ge 1\;\forall\;b>0\;\exists\;c>0\;\exists\;A\ge 1\;\forall\;\alpha\in\NN^d:\;\;\;M^{\frac{1}{b}}_{\alpha}\le AN^{\frac{1}{c}}\star L^{\frac{1}{a}}_{\alpha},
\end{equation}

\begin{equation}\label{PDOmainrelationgen1}
\forall\;a,b\ge 1\;\exists\;c\ge 1\;\exists\;A\ge 1\;\forall\;\alpha\in\NN^d:\;\;\;M^{\frac{1}{c}}_{\alpha}\le AN^{\frac{1}{b}}\star L^{\frac{1}{a}}_{\alpha}.
\end{equation}

\emph{Note:} Obviously, \eqref{PDOmainrelationgen} is equivalent to
\begin{equation}\label{PDOmainrelationgen2}
\forall\;a\ge 1\;\exists\;c>0\;\exists\;A\ge 1\;\forall\;\alpha\in\NN^d:\;\;\;M^{a}_{\alpha}\le AN^{\frac{1}{c}}\star L^{\frac{1}{a}}_{\alpha},
\end{equation}
and \eqref{PDOmainrelationgen1} is equivalent to
\begin{equation}\label{PDOmainrelationgen3}
\forall\;a\ge 1\;\exists\;c\ge 1\;\exists\;A\ge 1\;\forall\;\alpha\in\NN^d:\;\;\;M^{\frac{1}{c}}_{\alpha}\le AN^{\frac{1}{a}}\star L^{\frac{1}{a}}_{\alpha}.
\end{equation}

\begin{lemma}\label{PDOmainrelationgenlemma}
Let $\mathbf{M},\mathbf{N},\mathbf{L}\in\RR_{>0}^{\NN^d}$ be given weight sequences (cf. Definition \ref{weightsequdef}). Consider the following assertions:
\begin{itemize}
\item[$(i)$] The associated weight functions satisfy
\begin{equation}\label{PDOmainrelationgenasso}
\forall\;a\ge 1\;\forall\;b>0\;\exists\;c>0\;\exists\;B\ge 1\;\forall\;t\in\RR^d:\;\;\;\omega_{\mathbf{N}}(ct)+\omega_{\mathbf{L}}(at)\le\omega_{\mathbf{M}}(bt)+B.
\end{equation}

\item[$(ii)$] The associated weight functions satisfy
\begin{equation}\label{PDOmainrelationgenasso1}
\forall\;a\ge 1\;\forall\;b\ge 1\;\exists\;c\ge 1\;\exists\;B\ge 1\;\forall\;t\in\RR^d:\;\;\;\omega_{\mathbf{N}}(bt)+\omega_{\mathbf{L}}(at)\le\omega_{\mathbf{M}}(ct)+B.
\end{equation}
\end{itemize}
Then \eqref{PDOmainrelationgen} implies \eqref{PDOmainrelationgenasso} and \eqref{PDOmainrelationgen1} implies \eqref{PDOmainrelationgenasso1}. If $\mathbf{M}$ is in addition \emph{log-convex,} then \eqref{PDOmainrelationgenasso} implies \eqref{PDOmainrelationgen} and \eqref{PDOmainrelationgenasso1} implies \eqref{PDOmainrelationgen1}.
\end{lemma}

The proof shows that the parameters $a,b,c$ are preserved in all implications.

\demo{Proof}
Relation \eqref{PDOmainrelationgen} gives $\frac{N^{\frac{1}{c}}\star L^{\frac{1}{a}}_0|t^{\alpha}|}{N^{\frac{1}{c}}\star L^{\frac{1}{a}}_{\alpha}}\le\frac{N^{\frac{1}{c}}\star L^{\frac{1}{a}}_0}{M^{\frac{1}{b}}_0}A\frac{M^{\frac{1}{b}}_0|t^{\alpha}|}{M^{\frac{1}{b}}_{\alpha}}$ for all $t\in\RR^d$ and $\alpha\in\NN^d$ (the parameters $a,b,c$ are subject to \eqref{PDOmainrelationgen}). Thus the definition of associated weight functions, \eqref{anisoassofctadditivelemmaequ} and \eqref{217entire} together immediately give the conclusion with $B:=\log\frac{N^{\frac{1}{c}}\star L^{\frac{1}{a}}_0}{M^{\frac{1}{b}}_0}A=\log\frac{N_0L_0}{M_0}A$ and the same parameters $a,b,c$. When assuming \eqref{PDOmainrelationgen1}, then the argument is similar.\vspace{6pt}

On the other hand, \eqref{PDOmainrelationgenasso} combined with \eqref{assoweightfctinv}, \eqref{217entire}, \eqref{anisoassofctadditivelemmaequ} gives for all $\alpha\in\NN^d$:
\begin{align*}
M^{\frac{1}{b}}_{\alpha}&=M^{\frac{1}{b}}_0\sup_{t\in(0,+\infty)^d}\frac{t^{\alpha}}{\exp(\omega_{\mathbf{M}^{\frac{1}{b}}}(t))}=M_0\sup_{t\in(0,+\infty)^d}\frac{t^{\alpha}}{\exp(\omega_{\mathbf{M}}(bt))}
\\&
\le M_0e^B\sup_{t\in(0,+\infty)^d}\frac{t^{\alpha}}{\exp(\omega_{\mathbf{N}}(ct)+\omega_{\mathbf{L}}(at))}=M_0e^B\sup_{t\in(0,+\infty)^d}\frac{t^{\alpha}}{\exp(\omega_{\mathbf{N}^{\frac{1}{c}}}(t)+\omega_{\mathbf{L}^{\frac{1}{a}}}(t))}
\\&
=M_0e^B\sup_{t\in(0,+\infty)^d}\frac{t^{\alpha}}{\exp(\omega_{\mathbf{N}^{\frac{1}{c}}\star\mathbf{L}^{\frac{1}{a}}}(t))}=\frac{M_0e^B}{N_0L_0}N_0L_0\sup_{t\in(0,+\infty)^d}\frac{t^{\alpha}}{\exp(\omega_{\mathbf{N}^{\frac{1}{c}}\star\mathbf{L}^{\frac{1}{a}}}(t))}
\\&
=\frac{M_0e^B}{N_0L_0}(N^{\frac{1}{c}}\star L^{\frac{1}{a}})^{\on{lc}}_{\alpha}\le\frac{M_0e^B}{N_0L_0}N^{\frac{1}{c}}\star L^{\frac{1}{a}}_{\alpha}.
\end{align*}
Thus \eqref{PDOmainrelationgen} follows with the same parameters $a$, $b$, $c$ and $A:=\frac{M_0e^B}{N_0L_0}$. The implication \eqref{PDOmainrelationgenasso1}$\Rightarrow$\eqref{PDOmainrelationgen1} is analogous.
\qed\enddemo

As it is seen in Theorem \ref{Thm12igeneral}, the above conditions are optimal but they are quite involved. Therefore, we prove now that two ``easier-to-check'' resp. ``closed'' properties are sufficient.

\begin{lemma}\label{PDOmainrelationgenlemma1}
Let $\mathbf{M},\mathbf{N},\mathbf{L}\in\RR_{>0}^{\NN^d}$ be given.
\begin{itemize}
\item[$(i)$] The growth relation
\begin{equation}\label{PDOmainrelationgenvar}
\mathbf{M}\hyperlink{mtriangle}{\vartriangleleft}\mathbf{N}\star\mathbf{L}
\end{equation}
implies \eqref{PDOmainrelationgen} and \eqref{PDOmainrelationgen1}.

\item[$(ii)$] Assume that $\mathbf{M},\mathbf{N},\mathbf{L}$ are weight sequences (cf. Definition \ref{weightsequdef}) and such that
\begin{equation}\label{PDOmainrelationgenvar1}
(\mathbf{M},\mathbf{N})_{\on{mg}},\hspace{15pt}\mathbf{N}\hyperlink{mtriangle}{\vartriangleleft}\mathbf{L}.
\end{equation}
Then \eqref{PDOmainrelationgenasso} and \eqref{PDOmainrelationgenasso1} follow.
\end{itemize}
\end{lemma}

\demo{Proof}
$(i)$ Concerning \eqref{PDOmainrelationgen}, let $a\ge 1$ and $b>0$ be arbitrary, set $b_1:=\frac{b}{a}$ and so \eqref{PDOmainrelationgenvar} gives the existence of some $C_{b_1}=C_{a,b}\ge 1$ such that
$$\mathbf{M}_{\alpha}\le C_{b_1}b_1^{|\alpha|}\mathbf{N}\star\mathbf{L}_{\alpha}=C_{b_1}\min_{0\le\beta\le\alpha}\left(\frac{b}{a}\right)^{|\beta|}N_{\beta}\left(\frac{b}{a}\right)^{|\alpha|-|\beta|}L_{\alpha-\beta}=C_{b_1}b^{|\alpha|}N^{\frac{1}{a}}\star L^{\frac{1}{a}}_{\alpha}$$
for all $\alpha\in\NN^d$ which verifies \eqref{PDOmainrelationgen} with $c:=a$ and $A:=C_{b_1}$.

For \eqref{PDOmainrelationgen1}, we verify the equivalent condition \eqref{PDOmainrelationgen3}: Let $a\ge 1$ and put $b_1:=\frac{1}{a}$. Then \eqref{PDOmainrelationgen3} follows with $A:=C_{b_1}(=C_{a})$ and $c:=1$; i.e. $c$ can even be chosen uniformly for all $a$.\vspace{6pt}

$(ii)$ First, $(\mathbf{M},\mathbf{N})_{\on{mg}}$ precisely means that $\mathbf{M}\hyperlink{preceq}{\preceq}\mathbf{N}\star\mathbf{N}$, see $(d)$ in Section \ref{convolvesection}, and by $(c)$ there and the second part in \eqref{PDOmainrelationgenvar1} we obtain $\mathbf{M}\hyperlink{mtriangle}{\vartriangleleft}\mathbf{L}$, too. Let $a\ge 1$ and $b>0$ be arbitrary and $H\ge 1$ shall denote the constant appearing in \eqref{om6aniso} via $(\mathbf{M},\mathbf{N})_{\on{mg}}$; see Corollary \ref{anisoassofctadditivecor}. Then set $h:=\frac{b}{Ha}$ and via $\mathbf{N}\hyperlink{mtriangle}{\vartriangleleft}\mathbf{L}$ we involve \eqref{inclusioncharacterizationlemmaequ1} in Lemma \ref{inclusioncharacterizationlemma} in order to get the following estimate for all $t\in\RR^d$:
\begin{align*}
\omega_{\mathbf{N}}(aht)+\omega_{\mathbf{L}}(at)&\le\omega_{\mathbf{N}}(aht)+\omega_{\mathbf{N}}(aht)+C_h\le\omega_{\mathbf{M}}(Haht)+H+C_h
\\&
=\omega_{\mathbf{M}}(bt)+H+C_h.
\end{align*}
Thus \eqref{PDOmainrelationgenasso} follows with $c:=ah=\frac{b}{H}$ and $B:=H+C_h$ and note that $B$ is depending on both $a$ and $b$ (via $C_h$) whereas $c$ only depends on given $b$.

On the other hand, let $a,b\ge 1$, then put $h:=\frac{b}{a}$ in \eqref{inclusioncharacterizationlemmaequ1} and get for all $t\in\RR^d$:
\begin{align*}
\omega_{\mathbf{N}}(bt)+\omega_{\mathbf{L}}(at)&=\omega_{\mathbf{N}}(bt)+\omega_{\mathbf{L}}(bh^{-1}t)\le\omega_{\mathbf{N}}(bt)+\omega_{\mathbf{N}}(bt)+C_h
\\&
\le\omega_{\mathbf{M}}(Hbt)+H+C_h.
\end{align*}
Then \eqref{PDOmainrelationgenasso1} is shown with $B:=H+C_h$ and $c:=Hb$ (again only depending on $b$).
\qed\enddemo

\emph{Note:}

\begin{itemize}
\item[$(a)$] In view of $(c)$ in Section \ref{convolvesection} we have that \eqref{PDOmainrelationgenvar} for $\mathbf{M}=\mathbf{N}$ and/or $\mathbf{M}=\mathbf{L}$ has to fail since this would imply $\mathbf{M}\hyperlink{mtriangle}{\vartriangleleft}\mathbf{M}$ which is a contradiction for any $\mathbf{M}\in\RR_{>0}^{\NN^d}$. On the other hand, the relative growth between $\mathbf{N}$ and $\mathbf{L}$ in \eqref{PDOmainrelationgenvar} can behave irregular; more precisely it is not excluded to have even the strongest possible oscillation:
    $$0=\inf_{\alpha\in\NN^d\backslash\{0\}}\left(\frac{N_p}{L_p}\right)^{1/|\alpha|}<\sup_{\alpha\in\NN^d\backslash\{0\}}\left(\frac{N_p}{L_p}\right)^{1/|\alpha|}=+\infty.$$
    For an isotropic (counter-)example satisfying this behavior we refer to \cite[Sect. 3]{GelfandShilovincl}.

     Condition \eqref{PDOmainrelationgenvar} is tailor-made when dealing with the convolved matrices \eqref{convolvedmatrix} and \eqref{convolvedconstmatrix}.

\item[$(b)$] In the isotropic resp. one-dimensional situation and since (each) $\mathbf{G}^s$ satisfies moderate growth, i.e. $(\mathbf{G}^s,\mathbf{G}^s)_{\on{mg}}$, we have that \eqref{PDOmainrelationgenvar1} holds for $\mathbf{M}=\mathbf{G}^s=\mathbf{N}$ and $\mathbf{L}=\mathbf{G}^{s'}$ with $1<s<s'$. This corresponds to the particular case studied in \cite[Thm. 1.2 $(i)$]{OliaroPopivanov06}.

\item[$(c)$] More generally, \eqref{PDOmainrelationgenvar1} is easy-to-check when $\mathbf{M}$ satisfies moderate growth since then the first assumption in this condition is clear with $\mathbf{N}=\mathbf{M}$ and it only remains to find $\mathbf{L}$ such that the second part holds.
\end{itemize}

\subsection{Paley-Wiener-type result}\label{PWtheoremsection}
In order to follow the proof of \cite[Thm. 1.2 $(i)$]{OliaroPopivanov06} we have to recall the following crucial characterization \emph{(Paley-Wiener-type result);} see \cite[Lemma 3.3]{Komatsu73}.

\begin{lemma}\label{Komatsulemma33}
Let $\mathbf{M}=(M_{p})_{p\in\NN}$ be a log-convex weight sequence. Then the following hold:
\begin{itemize}
\item[$(i)$] Let $h>0$ and $f\in\mathcal{D}_{\mathbf{M},h}(K)$ with
$$\mathcal{D}_{\mathbf{M},h}(K):=\left\{f\in\mathcal{E}(\RR^d): \supp(f)\subseteq K,\;\;\|f\|_{\mathbf{M},K,h}:=\sup_{x\in K,\alpha\in\NN^d}\frac{|\partial^{\alpha}f(x)|}{h^{|\alpha|}M_{|\alpha|}}<+\infty\right\}.$$
Then the \emph{Fourier-Laplace transform} $\mathcal{F}(f)(\zeta):=\int_{\RR^d}f(x)\exp(-i\langle x,\zeta\rangle)dx$ is an entire function (on $\CC^d$) which satisfies on $\CC^d$ the estimate
    $$|\mathcal{F}(f)(\zeta)|\le M_0|K|\|f\|_{\mathbf{M},K,h}\exp(-\omega_{\mathbf{M}}(|\zeta|/(\sqrt{d}h))+H_K(\zeta)),$$
with $|K|$ denoting the \emph{Lebesgue measure} of $K$ and $H_K(\zeta):=\sup_{z\in K}\Im\langle z,\zeta\rangle$ is the \emph{support function} of $K$.

\item[$(ii)$] Let $h>0$ and assume that $\mathcal{F}(f)$ is a measurable function on $\RR^d$ such that $$\|\xi\mapsto\exp(\omega_{\mathbf{M}}(\xi/h))\mathcal{F}(f)(\xi)\|_{L^1(\RR^d)}<+\infty,$$
    then the \emph{inverse Fourier transform} $\mathcal{F}^{-1}\mathcal{F}(f)$ given by $$\mathcal{F}^{-1}(g)(x)=(2\pi)^{-d}\int_{\RR^d}g(\xi)\exp(i\langle\xi,x\rangle)d\xi$$
    belongs to the \emph{global ultradifferentiable class} $\mathcal{B}_{\mathbf{M},h}(\RR^d)$; i.e. $\mathcal{F}^{-1}(g)$ satisfies global estimates on whole $\RR^d$ meaning that $\|f\|_{\mathbf{M},\RR^d,h}<+\infty$.
\end{itemize}
\end{lemma}

Indeed, in \cite[Lemma 3.3]{Komatsu73} log-convexity was not assumed and in this case in $(ii)$ one has $\mathcal{B}_{\mathbf{M}^{\on{lc}},h}(\RR^d)$. And in $(i)$, when considering the Fourier transform instead and hence restricting to $\zeta\in\RR^d$, then $H_K(\zeta)=0$ and the estimate simplifies. Finally, note that for any weight matrix $\mathcal{M}=\{\mathbf{M}^{(\iota)}: \iota\in\mathcal{I}\}$ one gets the following representations as l.c.v.s. (with $K\subseteq\RR^d$ being compact regular; cf. \cite[Def. 2.5]{Komatsu73}):
\begin{equation}\label{Roumieutestfctrepr}
\mathcal{D}_{\{\mathcal{M}\}}(U)=\underset{K\subset\subset U}{\varinjlim}\;\underset{\iota\in\mathcal{I}}{\varinjlim}\;\underset{h>0}{\varinjlim}\;\mathcal{D}_{\mathbf{M}^{(\iota)},h}(K),
\end{equation}
\begin{equation}\label{Beurlingtestfctrepr}
\mathcal{D}_{(\mathcal{M})}(U)=\underset{K\subset\subset U}{\varinjlim}\;\underset{\iota\in\mathcal{I}}{\varprojlim}\;\underset{h>0}{\varprojlim}\;\mathcal{D}_{\mathbf{M}^{(\iota)},h}(K).
\end{equation}

\subsection{Main result}\label{hyperbolicmainthmsection}
The next definition generalizes \cite[Def. 1.1]{OliaroPopivanov06} to the mixed weight sequence setting and we allow for a loss of regularity expressed by a (different) sequence. This idea is natural when dealing with (associated) weight matrices due to the additional flexibility w.r.t. the index $\iota\in\mathcal{I}$.

\begin{definition}\label{locsovdef}
Let $U\subseteq\RR^d$ be open, $x_0\in U$, and let $\mathbf{M},\mathbf{N}\in\RR_{>0}^{\NN}$ be log-convex and non-quasianalytic weight sequences (cf. Definition \ref{weightsequdef}). A partial differential operator $P$ with coefficients belonging to the class $\mathcal{E}_{[\mathbf{M}]}(U)$ is called \emph{$[\mathbf{M},\mathbf{N}]$-locally solvable at $x_0\in U$} if there exists a neighborhood $V\subseteq U$ of $x_0$ such that for every $f\in\mathcal{D}_{[\mathbf{M}]}(V)$ there is a solution $u\in\mathcal{D}'_{[\mathbf{N}]}(V)$ of the equation $Pu=f$.
\end{definition}

Thus \cite[Def. 1.1]{OliaroPopivanov06} corresponds to the Roumieu-case $[\cdot,\cdot]=\{\cdot,\cdot\}$ and $\mathbf{M}=\mathbf{G}^s=\mathbf{N}$, $s>1$, and the assumptions on the sequences are natural since it follows that $\mathcal{D}_{[\mathbf{M}]}(V),\mathcal{D}_{[\mathbf{N}]}(V)\neq\{0\}$.

We formulate and prove now the main result of this section:

\begin{theorem}\label{Thm12igeneral}
Let $\mathbf{M}=(M_j)_{j\in\NN}$, $\mathbf{N}=(N_j)_{j\in\NN}$, $\mathbf{N^1}=(N^1_j)_{j\in\NN}$, $\mathbf{L}=(L_j)_{j\in\NN}$ be log-convex weight sequences and assume that $\mathbf{M}$ and $\mathbf{L}$ are non-quasianalytic and such that $(\mathbf{N},\mathbf{N}^1)_{\on{dc}}$, $(\mathbf{N}^1,\mathbf{L})_{\on{dc}}$. Let $p,q\in\NN$ and $k\in\NN_{>0}$ be given parameters such that $p>2q+1$, $p$ is even and $\frac{p-2q-1}{p-q}<\frac{1}{k}$. Set $s:=\frac{p-q}{k(p-2q-1)}(>1)$.
\begin{itemize}
\item[$(a)$] \emph{The Roumieu case.} Assume that
\begin{equation}\label{Thm12igeneralequ}
\forall\;a\ge 1\;\forall\;b>0\;\exists\;c>0\;\exists\;A\ge 1\;\forall\;j\in\NN:\;\;\;M^{\frac{1}{b}}_j\le AN^{\frac{1}{c}}\star (G^{s})^{\frac{1}{a}}_j;
\end{equation}
i.e. \eqref{PDOmainrelationgen} for the isotropic case and the special choice $\mathbf{L}=\mathbf{G}^s$ according to this parameter $s$. Then the operator
$$P(t,\partial_t,D_x):=\partial_t+at^pD_x^{2k}-bt^qD^k_x,\;\;\;a,b>0,$$
is $\{\mathbf{M},\mathbf{L}\}$-locally solvable at $(0,0)$.

\item[$(b)$] \emph{The Beurling case.} Assume that
\begin{equation}\label{Thm12igeneralequbeur}
\forall\;a\ge 1\;\forall\;b\ge 1\;\exists\;c\ge 1\;\exists\;A\ge 1\;\forall\;j\in\NN:\;\;\;M^{\frac{1}{c}}_j\le AN^{\frac{1}{b}}\star (G^{s})^{\frac{1}{a}}_j;
\end{equation}
i.e. \eqref{PDOmainrelationgen1} for the isotropic case and the special choice $\mathbf{L}=\mathbf{G}^s$, then the operator $P(t,\partial_t,D_x)$ is $(\mathbf{M},\mathbf{L})$-locally solvable at $(0,0)$.
\end{itemize}
\end{theorem}

\demo{Proof}
\emph{The Roumieu case.} Since $\mathbf{M}$ is log-convex, assumption \eqref{Thm12igeneralequ} is equivalent to
\begin{equation}\label{Thm12igeneralequ1}
\forall\;a\ge 1\;\forall\;b>0\;\exists\;c>0\;\exists\;C\ge 1\;\forall\;t\ge 0:\;\;\;\omega_{\mathbf{N}}(ct)+\omega_{\mathbf{G}^s}(at)\le\omega_{\mathbf{M}}(bt)+C;
\end{equation}
see Lemma \ref{PDOmainrelationgenlemma} (cf. \eqref{PDOmainrelationgenasso}). It is known that $\omega_{\mathbf{G}^s}\hyperlink{sim}{\sim}\id^{1/s}$ for any $s>0$, see e.g. \cite[Ex. 2.9]{genLegendreconj}, and hence we get that \eqref{Thm12igeneralequ1} implies
\begin{equation}\label{Thm12igeneralequ2}
\exists\;A\ge 1\;\forall\;a\ge 1\;\forall\;b>0\;\exists\;c>0\;\exists\;C\ge 1\;\forall\;t\ge 0:\;\;\;\omega_{\mathbf{N}}(ct)+A^{-1}(at)^{1/s}\le\omega_{\mathbf{M}}(bt)+C+A.
\end{equation}
The constant $A\ge 1$ is coming from the equivalence $\omega_{\mathbf{G}^s}\hyperlink{sim}{\sim}\id^{1/s}$ and depends only on the given and fixed parameter $s$.\vspace{6pt}

Now we follow the ideas in the proof of \cite[Thm. 1.2 $(i)$, p. 137-140]{OliaroPopivanov06}. We apply the \emph{partial Fourier transform} $\mathcal{F}_{x\mapsto\xi}$ to the equation
\begin{equation}\label{mainPDEequ}
P(t,\partial_t,D_x)u=f,
\end{equation}
and get (see \cite[$(2.2)$]{OliaroPopivanov06})
\begin{equation}\label{mainPDEequ1}
\partial_t\widehat{u}(\xi,t)+(at^p\xi^{2k}-bt^q\xi^k)\widehat{u}(\xi,t)=\widehat{f}(\xi,t),
\end{equation}
with $\widehat{u}(\xi,t)=\mathcal{F}_{x\mapsto\xi}[u(x,t)]$ and similarly for $\widehat{f}$. The solution of \eqref{mainPDEequ1} is then given as follows, see \cite[$(2.3)$]{OliaroPopivanov06}:
\begin{equation}\label{Oliaroequ23}
\widehat{u}(\xi,t)=\int_{-T}^t\exp(B(\xi,t,t'))\widehat{f}(\xi,t')dt',\;\;\;T>0\;\text{fixed.}
\end{equation}
Here, as in \cite{OliaroPopivanov06}, we have set
$$B(\xi,t,t'):=-\int_{t'}^t(ax^p\xi^{2k}-bx^q\xi^k)dx.$$
Finally, when applying the inverse partial Fourier transform to \eqref{Oliaroequ23} the formal solution of \eqref{mainPDEequ} is obtained by the following identity (see \cite[$(2.4)$]{OliaroPopivanov06}):
\begin{equation}\label{Oliaroequ24}
u(x,t)=\frac{1}{2\pi}\int_{-\infty}^{+\infty}e^{ix\xi}\widehat{u}(\xi,t)d\xi.
\end{equation}
Therefore, we have to estimate $\widehat{u}$ via \eqref{Oliaroequ23} and involve Lemma \ref{Komatsulemma33}. First, note that in \cite[Lemma 2.1]{OliaroPopivanov06} it is shown that
\begin{equation}\label{Oliaroequ25}
\exists\;C_1>0\;\forall\;\xi\;\forall\;t'\in[-T,t]:\;\;\;B(\xi,t,t')\le C_1|\xi|^{1/s}.
\end{equation}
We use the estimate in the proof of \cite[Thm. 1.2 $(i)$, p. 140]{OliaroPopivanov06} and $(i)$ in Lemma \ref{Komatsulemma33} for $f$, i.e. involving the Roumieu-type ultradifferentiable regularity w.r.t. the sequence $\mathbf{M}$ since by assumption $f\in\mathcal{D}_{\{\mathbf{M}\}}(V)$ and $V\subseteq\RR^2$ being a neighborhood of $(0,0)$. Thus, by recalling also the representation \eqref{Roumieutestfctrepr}, there exist $c>0$ and $D>0$ such that
$$|\widehat{u}(\xi,t)|\le D\int_{-T}^T\exp\left(C_1|\xi|^{1/s}\right)\exp(-\omega_{\mathbf{M}}(c|\xi|))dt'\le 2DT\exp\left(C_1|\xi|^{1/s}-\omega_{\mathbf{M}}(c|\xi|)\right).$$
Here the constant $C_1$ is coming from \eqref{Oliaroequ25}, whereas $D$ and $c$ are subject to $(i)$ in Lemma \ref{Komatsulemma33} and note that we are dealing with a real variable and hence the additional term in the exponent with the supporting function $H_K$ disappears. We apply \eqref{Thm12igeneralequ2} to the choices $b:=c$ and $a:=(C_1A)^s$, which is possible since $A$ is only depending on $s$, and then for some constant $c_1>0$ being related to $b$ and $a$ via \eqref{Thm12igeneralequ2}:
\begin{equation}\label{mainthmfinalequ}
2DT\exp\left(C_1|\xi|^{1/s}-\omega_{\mathbf{M}}(c|\xi|)\right)\le 2DT\exp\left(-\omega_{\mathbf{N}}(c_1|\xi|)+C+A\right).
\end{equation}
This estimate and the known fact that $\log(t)=o(\omega_{\mathbf{N}}(t))$ as $t\rightarrow+\infty$ (i.e. \hyperlink{om3}{$(\omega_3)$}) imply that $u$ given by \eqref{Oliaroequ24} is well defined: indeed the integrand $\widehat{u}(\xi,t)$ in \eqref{Oliaroequ24} is exponentially decreasing as $|\xi|\rightarrow+\infty$ and thus the integral is finite. For this recall that by definition for any weight sequence $\mathbf{N}$ we do have \hyperlink{om3}{$(\omega_3)$} for $\omega_{\mathbf{N}}$ and, moreover, for any arbitrary (continuous) function $\omega:[0,+\infty)\rightarrow[0,+\infty)$ satisfying \hyperlink{om3}{$(\omega_3)$} one has that for each $n\in\NN_{>0}$ there exists $t_n>0$ such that $\log(t)\le\frac{1}{n}\omega(t)$. This gives $\frac{1}{t^n}\ge\exp(-\omega(t))$ for all $t\ge t_n$ and hence finiteness of the integral.

Finally, by relation $(\mathbf{N},\mathbf{N}^1)_{\on{dc}}$ we obtain
\begin{equation}\label{mixeddcassofct}
\exists\;D_1\ge 1\;\exists\;D_2>0\;\forall\;t>0:\;\;\;\omega_{\mathbf{N^1}}(t)+\log(t)\le\omega_{\mathbf{N}}(D_1t)+D_2,
\end{equation}
and the analogous estimate holds via $(\mathbf{N}^1,\mathbf{L})_{\on{dc}}$. We combine both estimates and so
\begin{equation}\label{mixeddcassofctcombine}
\exists\;D_1\ge 1\;\exists\;D_2>0\;\forall\;t>0:\;\;\;\omega_{\mathbf{L}}(t)+2\log(t)\le\omega_{\mathbf{N}}(D_1t)+D_2.
\end{equation}
The proof that $(\mathbf{N},\mathbf{N}^1)_{\on{dc}}$ is equivalent to \eqref{mixeddcassofct} can be found in \cite[Lemma 2]{nuclearglobal2} (see also \cite[Lemma 3.4]{Komatsu73}) and this step should be compared with the proof of \cite[Prop. 1]{nuclearglobal2} (for the case $d=1$). Hence, by combining \eqref{mainthmfinalequ} and \eqref{mixeddcassofctcombine} it follows that
$$\exp\left(\omega_{\mathbf{L}}\left(\frac{c_1|\xi|}{D_1}\right)\right)|\widehat{u}(\xi,t)|\le 2DT\exp(A+B+D_2)\exp\left(-2\log\left(\frac{c_1|\xi|}{D_1}\right)\right).$$
And this estimate implies that the assumption in $(ii)$ Lemma \ref{Komatsulemma33} for the partial Fourier transform $\mathcal{F}_{x\mapsto\xi}(u)$ holds w.r.t. $\mathbf{L}$ and $h:=\frac{D_1}{c_1}$. Thus $(ii)$ Lemma \ref{Komatsulemma33} yields the conclusion in the Roumieu case. Note that this final estimate holds uniformly for any $t$.\vspace{6pt}

\emph{The Beurling case.} Similarly, \eqref{Thm12igeneralequbeur} implies
\begin{equation}\label{Thm12igeneralequ2beur}
\exists\;A\ge 1\;\forall\;a\ge 1\;\forall\;b\ge 1\;\exists\;c\ge 1\;\exists\;C\ge 1\;\forall\;t\ge 0:\;\;\;\omega_{\mathbf{N}}(bt)+A^{-1}(at)^{1/s}\le\omega_{\mathbf{M}}(ct)+C+A.
\end{equation}
Then follow the lines in the Roumieu case and note that now by $(i)$ in Lemma \ref{Komatsulemma33} (and by taking into account \eqref{Beurlingtestfctrepr}) even for any $c>0$ we can find $D>0$ such that
\begin{equation}\label{Thm12igeneralequ2beur1}
|\widehat{u}(\xi,t)|\le D\int_{-T}^T\exp\left(C_1|\xi|^{1/s}\right)\exp(-\omega_{\mathbf{M}}(c|\xi|))dt'\le 2DT\exp\left(C_1|\xi|^{1/s}-\omega_{\mathbf{M}}(c|\xi|)\right).
\end{equation}
In order to conclude, let $h>0$ be arbitrary (but fixed), set $c_1:=\frac{D_1}{h}$ and then via \eqref{Thm12igeneralequ2beur} applied to $b:=c_1$ and $a:=(C_1A)^s$ one can find $c>0$ such that \eqref{mainthmfinalequ} is valid and \eqref{Thm12igeneralequ2beur1} holds with this given $c>0$ and some $D=D(c)>0$. Summarizing, this shows via $(ii)$ in Lemma \ref{Komatsulemma33} the precise Beurling-type regularity of the solution $u$ in terms of $\mathbf{L}$.
\qed\enddemo

An immediate consequence of the previous statement is \cite[Thm. 1.2 $(i)$]{OliaroPopivanov06} and we see that this result can be transferred to the Beurling setting, too.

\begin{corollary}
Let $p,q\in\NN$, $k\in\NN_{>0}$, $s:=\frac{p-q}{k(p-2q-1)}(>1)$ and the operator $P(t,\partial_t,D_x)$ be as in Theorem \ref{Thm12igeneral}. Then this operator is $[\mathbf{G}^{s'},\mathbf{G}^{s'}]$-locally solvable at $(0,0)$ for any $1<s'<s$.
\end{corollary}

\demo{Proof}
Since each $\mathbf{G}^s$ satisfies moderate growth (and hence derivation closedness too) we apply Theorem \ref{Thm12igeneral} to $\mathbf{M}=\mathbf{N}=\mathbf{N}^1=\mathbf{L}=\mathbf{G}^{s'}$ and $(ii)$ in Lemma \ref{PDOmainrelationgenlemma1} to $\mathbf{M}=\mathbf{N}=\mathbf{G}^{s'}$ and $\mathbf{L}=\mathbf{G}^s$; see also comment $(b)$ in Section \ref{preliminarysection}. Note that Lemma \ref{PDOmainrelationgenlemma1} implies both the Roumieu- and the Beurling-type estimates and hence this result is valid for the Beurling setting, too.
\qed\enddemo

In view of Lemmas \ref{PDOmainrelationgenlemma} and \ref{PDOmainrelationgenlemma1}, both \eqref{PDOmainrelationgenvar} and \eqref{PDOmainrelationgenvar1} have immediate consequences and generalizations of Theorem \ref{Thm12igeneral} in the weight matrix setting; i.e. when involving a matrix parameter $\iota$. More precisely, for a fixed sequence $\mathbf{L}$ and weight matrices $\mathcal{M}:=\{\mathbf{M}^{(\iota)}: \iota\in\mathcal{I}\}$, $\mathcal{N}:=\{\mathbf{N}^{(\iota)}: \iota\in\mathcal{I}\}$ let us consider
\begin{equation}\label{PDOmainrelationgenvarRoum}
\forall\;\iota\in\mathcal{I}\;\exists\;\iota_1\in\mathcal{I}:\;\;\;\mathbf{M}^{(\iota)}\hyperlink{mtriangle}{\vartriangleleft}\mathbf{N}^{(\iota_1)}\star\mathbf{L},
\end{equation}
\begin{equation}\label{PDOmainrelationgenvarBeur}
\forall\;\iota\in\mathcal{I}\;\exists\;\iota_1\in\mathcal{I}:\;\;\;\mathbf{M}^{(\iota_1)}\hyperlink{mtriangle}{\vartriangleleft}\mathbf{N}^{(\iota)}\star\mathbf{L},
\end{equation}
and recognize that on the right-hand sides the convolved sequence is corresponding to the elements of the matrix $\mathcal{N}\star\mathbf{L}$; see \eqref{convolvedconstmatrix}.

\begin{theorem}\label{Thm12igeneral1}
Let $\mathcal{M}:=\{\mathbf{M}^{(\iota)}: \iota\in\mathcal{I}\}$, $\mathcal{L}:=\{\mathbf{L}^{(\iota)}: \iota\in\mathcal{I}\}$ be standard log-convex and non-quasianalytic weight matrices, and let $\mathcal{N}:=\{\mathbf{N}^{(\iota)}: \iota\in\mathcal{I}\}$, $\mathcal{Q}:=\{\mathbf{Q}^{(\iota)}: \iota\in\mathcal{I}\}$ be standard log-convex weight matrices. Let $p,q\in\NN$, $k\in\NN_{>0}$ and $s:=\frac{p-q}{k(p-2q-1)}(>1)$ be as in Theorem \ref{Thm12igeneral}.
\begin{itemize}
\item[$(i)$] Assume that \eqref{PDOmainrelationgenvarRoum} holds with $\mathbf{L}=\mathbf{G}^s$ and for the matrices $\mathcal{M}$, $\mathcal{N}$ and that
$$\forall\;\iota,\iota_1\in\mathcal{I}\;\exists\;\iota_2,\iota_3\in\mathcal{I}:\;\;\;(\mathbf{N}^{(\iota)},\mathbf{Q}^{(\iota_2)})_{\on{dc}},\;(\mathbf{Q}^{(\iota_1)},\mathbf{L}^{(\iota_3)})_{\on{dc}}.$$

Then for all $\iota\in\mathcal{I}$ there exists $\iota_1\in\mathcal{I}$ such that the operator $P(t,\partial_t,D_x)$ is $[\mathbf{M}^{(\iota)},\mathbf{L}^{(\iota_1)}]$-locally solvable at $(0,0)$.

\item[$(ii)$] Assume that \eqref{PDOmainrelationgenvarBeur} holds with $\mathbf{L}=\mathbf{G}^s$ and for the matrices $\mathcal{M}$, $\mathcal{N}$ and that
$$\forall\;\iota,\iota_1\in\mathcal{I}\;\exists\;\iota_2,\iota_3\in\mathcal{I}:\;\;\;(\mathbf{N}^{(\iota_3)},\mathbf{Q}^{(\iota_1)})_{\on{dc}},\;(\mathbf{Q}^{(\iota_2)},\mathbf{L}^{(\iota)})_{\on{dc}}$$
Then for all $\iota\in\mathcal{I}$ there exists $\iota_1\in\mathcal{I}$ such that the operator $P(t,\partial_t,D_x)$ is $[\mathbf{M}^{(\iota_1)},\mathbf{L}^{(\iota)}]$-locally solvable at $(0,0)$.
\end{itemize}
\end{theorem}

\demo{Proof}
The conclusion follows immediately by (the proof of) Theorem \ref{Thm12igeneral}, Lemma \ref{PDOmainrelationgenlemma} and $(i)$ in Lemma \ref{PDOmainrelationgenlemma1} and recall that we have in both cases $[\cdot,\cdot]$-local solvability since $(i)$ in Lemma \ref{PDOmainrelationgenlemma1} implies both the Roumieu-type and the Beurling-type condition.
\qed\enddemo

Theorem \ref{Thm12igeneral1} simplifies when one can take $\mathcal{M}=\mathcal{Q}=\mathcal{L}$; it is desirable for applications to consider weight matrices which are closed under $(\cdot,\cdot)_{\on{dc}}$ (cf. \hyperlink{R-mg}{$(\mathcal{M}_{\{\on{mg}\}})$} and \hyperlink{B-mg}{$(\mathcal{M}_{(\on{mg})})$}). This is ensured, in particular, in the BMT-weight function setting and we focus now on this case and exploit the information from Corollary \ref{convolvedmatrixthmcor}.

\begin{theorem}\label{Thm12igeneral3}
Let $\sigma,\tau\in\hyperlink{omset0}{\mathcal{W}_0}$ be given with associated weight matrices $\mathcal{M}_{\sigma}:=\{\mathbf{S}^{(\ell)}: \ell>0\}$, $\mathcal{M}_{\tau}:=\{\mathbf{T}^{(\ell)}: \ell>0\}$ and assume that both weights satisfy \hyperlink{om1}{$(\omega_1)$}. Then we get:
\begin{itemize}
\item[$(i)$] For any $\alpha>0$, the following equivalence is valid:
\begin{equation}\label{Thm12igeneral3equ}
\mathcal{M}_{\sigma}\vartriangleleft\mathcal{M}_{\tau}\star\mathbf{G}^{\alpha}\Leftrightarrow\tau(t)+t^{1/\alpha}=o(\sigma(t)),\;t\rightarrow+\infty;
\end{equation}
i.e. $\mathcal{M}_{\sigma}\vartriangleleft\mathcal{M}_{\tau}\star\mathbf{G}^{\alpha}$ if and only if $\sigma\hyperlink{omvartriangle}{\vartriangleleft}\tau+\id^{1/\alpha}$.

\item[$(ii)$] Let $p,q\in\NN$, $k\in\NN_{>0}$ and $s:=\frac{p-q}{k(p-2q-1)}(>1)$ be as in Theorem \ref{Thm12igeneral}. Assume that $\sigma$ and $\tau$ are also non-quasianalytic and that \eqref{Thm12igeneral3equ} is valid for $s=\alpha$. Then, for any $\ell,\ell_1>0$, the operator $P(t,\partial_t,D_x)$ is $[\mathbf{S}^{(\ell)},\mathbf{T}^{(\ell_1)}]$-locally solvable at $(0,0)$.
\end{itemize}
\end{theorem}

\demo{Proof}
$(i)$ First, by the aforementioned equivalence $\omega_{\mathbf{G}^{\alpha}}\hyperlink{sim}{\sim}\id^{1/\alpha}$ we have that $\tau(t)+\omega_{\mathbf{G}^{\alpha}}(t)=o(\sigma(t))$ if and only if $\tau(t)+t^{1/\alpha}=o(\sigma(t))$. This equivalence also yields that (each) $\omega_{\mathbf{G}^{\alpha}}$ satisfies \hyperlink{om1}{$(\omega_1)$} and \hyperlink{om6}{$(\omega_6)$}. Of course, this can be checked directly since (each) $\mathbf{G}^{\alpha}$ satisfies moderate growth, see e.g. \cite[Thm. 3.1]{modgrowthstrange} and the references there, and concerning \hyperlink{om1}{$(\omega_1)$} one can check that \cite[Thm. 3.1]{subaddlike} can be applied to $\mathbf{G}^{\alpha}$ since \cite[$(3.2)$]{subaddlike} holds with any $L\in\NN_{\ge 3}$ (by the estimates $p!\le p^p\le e^p$ for all $p\in\NN$ and convention $0^0:=1$).

Second, Lemma \ref{assoweightomega0} applied to $\mathbf{G}^{\alpha}\in\hyperlink{LCset}{\mathcal{LC}}$ gives $\mathbf{G}^{\alpha}=(\mathbf{G}^{\alpha})^{(1)}$ (recall again \cite[$(2.13)$]{subaddlike}) and that the matrix $\mathcal{M}_{\omega_{\mathbf{G}^{\alpha}}}$ is constant, hence $\mathcal{M}_{\tau}\star\mathcal{M}_{\omega_{\mathbf{G}^{\alpha}}}$ and $\mathcal{M}_{\tau}\star\mathbf{G}^{\alpha}$ are $R$- and $B$-equivalent which immediately gives (as l.c.v.s.) $\mathcal{E}_{[\mathcal{M}_{\tau}\star\mathcal{M}_{\omega_{\mathbf{G}^{\alpha}}}]}=\mathcal{E}_{[\mathcal{M}_{\tau}\star\mathbf{G}^{\alpha}]}$. And by Corollary \ref{convolvedmatrixthmcor} we get $\mathcal{E}_{[\tau+\omega_{\mathbf{G}^{\alpha}}]}=\mathcal{E}_{[\mathcal{M}_{\tau}\star\mathbf{G}^{\alpha}]}$. The same equalities hold for all analogously defined weighted spaces, in particular for the corresponding test function space which is non-trivial if and only if $\tau$ is non-quasianalytic and $\alpha>1$; see $(i)$ in Corollary \ref{convolvedmatrixthmcor1}.

Then $\mathcal{M}_{\sigma}\vartriangleleft\mathcal{M}_{\tau}\star\mathbf{G}^{\alpha}$ implies the (continuous) inclusion $\mathcal{E}_{\{\sigma\}}=\mathcal{E}_{\{\mathcal{M}_{\sigma}\}}\subseteq\mathcal{E}_{(\mathcal{M}_{\tau}\star\mathbf{G}^{\alpha})}=\mathcal{E}_{(\tau+\omega_{\mathbf{G}^{\alpha}})}$ and hence $\tau(t)+\omega_{\mathbf{G}^{\alpha}}(t)=o(\sigma(t))$ as $t\rightarrow+\infty$ follows by taking into account \cite[Cor. 5.17 (2)]{compositionpaper}. For this note that both weights $\sigma$ and $\tau+\omega_{\mathbf{G}^{\alpha}}$ belong to \hyperlink{omset0}{$\mathcal{W}_0$} and satisfy \hyperlink{om1}{$(\omega_1)$}.

Conversely, $\tau(t)+\omega_{\mathbf{G}^{\alpha}}(t)=o(\sigma(t))$ implies $\mathcal{M}_{\sigma}\vartriangleleft\mathcal{M}_{\tau+\omega_{\mathbf{G}^{\alpha}}}$ by \cite[Lemma 5.8 $(i)$]{genLegendreconjBMT} (indeed, each appearing weight in this relation satisfies \hyperlink{om1}{$(\omega_1)$}). Then, in view of \eqref{weightmatrixequivequ} the relation $\mathcal{M}_{\sigma}\vartriangleleft\mathcal{M}_{\tau}\star\mathcal{M}_{\omega_{\mathbf{G}^{\alpha}}}$ follows. Finally, since $\mathcal{M}_{\omega_{\mathbf{G}^{\alpha}}}$ is constant and $\mathbf{G}^{\alpha}=(\mathbf{G}^{\alpha})^{(1)}$ condition $\mathcal{M}_{\sigma}\vartriangleleft\mathcal{M}_{\tau}\star\mathbf{G}^{\alpha}$ is shown.

Note: When $\mathbf{G}^{\alpha}$ is replaced by any other sequence $\mathbf{N}\in\hyperlink{LCset}{\mathcal{LC}}$ satisfying moderate growth and \cite[$(3.2)$]{subaddlike}, then all arguments transfer immediately (except having the equivalence $\omega_{\mathbf{N}}\hyperlink{sim}{\sim}\id^{1/\alpha}$ for some $\alpha>0$ and hence replace $\id^{1/\alpha}$ by $\omega_{\mathbf{N}}$ in \eqref{Thm12igeneral3equ}).\vspace{6pt}

$(ii)$ This follows as Theorem \ref{Thm12igeneral1} with $\mathbf{L}=\mathbf{G}^s$, $\mathcal{M}=\mathcal{M}_{\sigma}$, $\mathcal{N}=\mathcal{Q}=\mathcal{L}=\mathcal{M}_{\tau}$: First, \eqref{newmoderategrowth} implies $(\mathbf{T}^{(\ell)},\mathbf{T}^{(2\ell)})_{\on{dc}}$ for all $\ell>0$ and, second, recall $\mathcal{M}_{\sigma}\vartriangleleft\mathcal{M}_{\tau}\star\mathbf{G}^s$ means
$$\forall\;\ell,\ell_1>0:\;\;\;\mathbf{S}^{(\ell)}\hyperlink{mtriangle}{\vartriangleleft}\mathbf{T}^{(\ell_1)}\star\mathbf{G}^s,$$
which is much stronger than \eqref{PDOmainrelationgenvarRoum} and \eqref{PDOmainrelationgenvarBeur} with $\mathbf{L}=\mathbf{G}^s$.
\qed\enddemo

\begin{theorem}\label{Thm12igeneral2}
Let $\omega\in\hyperlink{omset0}{\mathcal{W}_0}$ be given with associated weight matrix $\mathcal{M}_{\omega}:=\{\mathbf{W}^{(\ell)}: \ell>0\}$ and let $p,q\in\NN$, $k\in\NN_{>0}$ and $s:=\frac{p-q}{k(p-2q-1)}(>1)$ be as in Theorem \ref{Thm12igeneral}. Assume that $\omega$ is non-quasianalytic and that $t^{1/s}=o(\omega(t))$ as $t\rightarrow+\infty$. Then, for any $\ell>0$, the operator $P(t,\partial_t,D_x)$ is $[\mathbf{W}^{(\ell)},\mathbf{W}^{(8\ell)}]$-locally solvable at $(0,0)$.
\end{theorem}

\demo{Proof}
$\mathcal{M}_{\omega}$ is standard log-convex and non-quasianalytic and by $(ii)$ in Section \ref{weightfctsection} (see \eqref{newmoderategrowth}) we have $(\mathbf{W}^{(\ell)},\mathbf{W}^{(2\ell)})_{\on{mg}}$ and hence $(\mathbf{W}^{(\ell)},\mathbf{W}^{(2\ell)})_{\on{dc}}$ for all $\ell>0$. We apply Theorem \ref{Thm12igeneral} to $\mathbf{M}=\mathbf{W}^{(\ell)}$, $\mathbf{N}=\mathbf{W}^{(2\ell)}$, $\mathbf{N}^1=\mathbf{W}^{(4\ell)}$, $\mathbf{L}=\mathbf{W}^{(8\ell)}$ and $(ii)$ in Lemma \ref{PDOmainrelationgenlemma1} to $\mathbf{M}=\mathbf{W}^{(\ell)}$, $\mathbf{N}=\mathbf{W}^{(2\ell)}$, $\mathbf{L}=\mathbf{G}^s$ (see \eqref{PDOmainrelationgenvar1}). Note that this result yields both the Roumieu- and the Beurling-type requirement and thus we have $[\cdot,\cdot]$-local solvability. Therefore, it remains to treat the second part in \eqref{PDOmainrelationgenvar1} and we show $\mathbf{W}^{(\ell)}\hyperlink{mtriangle}{\vartriangleleft}\mathbf{G}^s$ for all $\ell>0$.\vspace{6pt}

Assumption $t^{1/s}=o(\omega(t))$, \eqref{goodequivalenceclassic} and the aforementioned equivalence $\omega_{\mathbf{G}^s}\hyperlink{sim}{\sim}\id^{1/s}$ imply $\omega_{\mathbf{G}^s}(t)=o(\omega_{\mathbf{W}^{(\ell)}}(t))$, i.e. $\omega_{\mathbf{W}^{(\ell)}}\hyperlink{omvartriangle}{\vartriangleleft}\omega_{\mathbf{G}^s}$, for all $\ell>0$. $\omega_{\mathbf{G}^s}$ satisfies \hyperlink{om1}{$(\omega_1)$} (cf. the proof of Theorem \ref{Thm12igeneral3}) and by \cite[Prop. 2.2 $(iii)$]{ultradifferentiablecomparison} we get relation $\omega_{\mathbf{W}^{(\ell)}}\vartriangleleft_{\mathfrak{c}}\omega_{\mathbf{G}^s}$ for each $\ell>0$; i.e.  \eqref{inclusioncharacterizationlemmaequ1} with $\mathbf{M}=\mathbf{W}^{(\ell)}$ and $\mathbf{N}=\mathbf{G}^s$. Then $(ii)$ in Lemma \ref{inclusioncharacterizationlemma} yields the conclusion.
\qed\enddemo

Note that the proof of Theorem \ref{Thm12igeneral2} does not require \hyperlink{om1}{$(\omega_1)$} for $\omega$ necessarily and hence the equality $\mathcal{D}_{[\mathcal{M}_{\omega}]}=\mathcal{D}_{[\omega]}$ is unclear in general. However, the previous result purely involves the sequences $\mathbf{W}^{(\ell)}$ and does not need this equality.

Non-quasianalyticity implies $\omega(t)=o(t)$, see e.g. \cite[Sect. 2.4 $(h)$]{ultradifferentiablecomparison}, and therefore $\omega$ in the previous result has to satisfy necessarily specific growth restrictions: $t^{1/s}=o(\omega(t))$ for a certain $s>1$ and $\omega(t)=o(t)$. We construct explicitly a \emph{non-constant} (associated) weight matrix to which Theorem \ref{Thm12igeneral2} can be applied.

\begin{proposition}\label{counterexample}
For any $s>1$ there exists a weight matrix $\mathcal{M}=\mathcal{M}_{\omega_{\mathbf{M}}}:=\{\mathbf{M}^{(\iota)}: \iota\in\mathcal{I}\}$ with the following properties:
\begin{itemize}
\item[$(a)$] $\mathcal{M}$ is standard log-convex,

\item[$(b)$] $\mathcal{M}$ is non-quasianalytic,

\item[$(c)$] $\mathcal{M}$ is non-constant,

\item[$(d)$] the following equivalent conditions hold:
\begin{itemize}
\item[$(*)$] $t^{1/s}=o(\omega_{\mathbf{M}}(t))$ as $t\rightarrow+\infty$,

\item[$(*)$] $\mathcal{M}\vartriangleleft\mathcal{M}_{\omega_{\mathbf{G}^s}}$,

\item[$(*)$] $\mathcal{M}\vartriangleleft\mathbf{G}^s$; i.e.
\begin{equation}\label{counterexampleequ}
	\forall\;\iota\in\mathcal{I}:\;\;\;\mathbf{M}^{(\iota)}\hyperlink{mtriangle}{\vartriangleleft}\mathbf{G}^s.
\end{equation}
\end{itemize}

\item[$(e)$] $\mathcal{M}$ satisfies \hyperlink{R-mg}{$(\mathcal{M}_{\{\on{mg}\}})$}, \hyperlink{B-mg}{$(\mathcal{M}_{(\on{mg})})$}.
\end{itemize}
Indeed, $\mathcal{M}$ is associated with a special weight function $\omega\in\hyperlink{omset0}{\mathcal{W}_0}$, namely with an associated weight function given in terms of $\mathbf{M}\in\RR_{>0}^{\NN}$. Thus $\omega=\omega_{\mathbf{M}}$ and this weight is as required in Theorem \ref{Thm12igeneral2}.
\end{proposition}

\demo{Proof}
Let $s>1$ be given and fixed. The idea is to define $\mathcal{M}$ in terms of $\mathcal{M}_{\omega_{\mathbf{M}}}=\{\mathbf{M}^{(\iota)}: \iota>0\}$ via $\mathbf{M}\in\hyperlink{LCset}{\mathcal{LC}}$ satisfying appropriate growth requirements (recall Lemma \ref{assoweightomega0}). In order to do so, we construct $\mathbf{M}$ via the corresponding sequence of quotients $\mu=(\mu_p)_{p\in\NN_{>0}}$; i.e. we set $M_p:=\prod_{i=1}^p\mu_i$, $p\in\NN$, and so $M_0:=1$ (empty product).\vspace{6pt}

Let $\mathbf{a}=(a_j)_{j\in\NN_{>0}}$ be a sequence in $\NN$ and let $\mathbf{b}=(b_j)_{j\in\NN_{>0}}$ be a sequence in $\RR_{>0}$ such that
\begin{itemize}
\item[$(i)$] $a_1=1=b_1$, both sequences are strictly increasing with $\lim_{j\rightarrow+\infty}b_j=+\infty$;

\item[$(ii)$] $j\mapsto\frac{a_j^s}{b_j}$ is non-decreasing and $\lim_{j\rightarrow+\infty}\frac{a_j^s}{b_j}=+\infty$;

\item[$(iii)$] $a_{j+1}\le\frac{a_j^s}{2^jb_j}+a_j$ for all $j\in\NN_{>0}$ sufficiently large;

\item[$(iv)$] $\lim_{j\rightarrow+\infty}\frac{a_{j+1}}{a_j}=+\infty$;

\item[$(v)$] $\sup_{j\in\NN_{>0}}\frac{a_{j+1}^s/b_{j+1}}{a_j^s/b_j}=+\infty$.
\end{itemize}

An explicit choice for $\mathbf{a}$ and $\mathbf{b}$ is given when taking, e.g., $a_j:=j!$ and $b_j:=j$, $j\in\NN_{>0}$: $(i)$ and $(iv)$ are clear, concerning $(ii)$ note that obviously $\frac{j!^s}{j}\rightarrow+\infty$ and $\frac{j!^s}{j}\le\frac{(j+1)!^s}{j+1}\Leftrightarrow 1+\frac{1}{j}\le(j+1)^s$ is clear for all $j\in\NN_{>0}$. And $(iii)$ follows since $(j+1)!\le\frac{j!^s}{2^jj}+j!\Leftrightarrow j+1\le\frac{j!^{s-1}}{2^jj}+1\Leftrightarrow 2^jj^2\le j!^{s-1}$ is valid for all $j$ large enough (note that $s>1$). Finally, $\frac{a_{j+1}^s/b_{j+1}}{a_j^s/b_j}=\frac{j(j+1)^s}{j+1}=j(j+1)^{s-1}\rightarrow+\infty$ and so $(v)$ holds. Note that these choices for $\mathbf{a}$ and $\mathbf{b}$ work \emph{uniformly} for \emph{any} $s>1$ and are not depending on this parameter. Then put
\begin{equation}\label{exampledef}
\mu_0:=1,\hspace{15pt}\mu_p:=\frac{a_j^s}{b_j},\hspace{5pt}a_j\le p\le a_{j+1}-1,\;\;\;j\in\NN_{>0}.
\end{equation}

$(a)$ By $(ii)$ the sequence $\mathbf{M}$ is log-convex with $\lim_{p\rightarrow+\infty}\mu_p=\lim_{p\rightarrow+\infty}(M_p)^{1/p}=+\infty$. Moreover, $M_0=1$ is clear by definition and $M_1=\mu_1=\mu_{a_1}=\frac{a_1^s}{b_1}=1$; see $(i)$. Hence $\mathbf{M}\in\hyperlink{LCset}{\mathcal{LC}}$ and so $\mathcal{M}_{\omega_{\mathbf{M}}}$ is standard log-convex; see Lemma \ref{assoweightomega0} and $(i)$ in Section \ref{weightfctsection}.\vspace{6pt}

$(b)$ First, we show that $\mathcal{M}_{\omega_{\mathbf{M}}}$ is non-quasianalytic if and only if $\mathbf{M}$ is non-quasianalytic and this should be compared with $(vi)$ in Section \ref{weightfctsection}. By \cite[Lemma 4.1]{Komatsu73} a sequence
$\mathbf{M}\in\hyperlink{LCset}{\mathcal{LC}}$ is non-quasianalytic if and only if $\omega_{\mathbf{M}}$ is non-quasianalytic. In view of \eqref{goodequivalenceclassic} and the fact that non-quasianalyticity is preserved under equivalence of weight functions, this holds if and only if some/each $\omega_{\mathbf{M}^{(\iota)}}$ is non-quasianalytic. Finally, \cite[Lemma 4.1]{Komatsu73} applied to some/each $\mathbf{M}^{(\iota)}$ yields that this holds if and only if some/each $\mathbf{M}^{(\iota)}$ is non-quasianalytic; i.e. $\mathcal{M}_{\omega_{\mathbf{M}}}$ is non-quasianalytic.

And $\mathbf{M}$ is non-quasianalytic since in view of $(iii)$ there exists $j_0\in\NN_{>0}$ such that the estimate there is valid for all $j>j_0$. And so, by \eqref{exampledef}, one has
\begin{align*}
\sum_{p\ge 1}\frac{1}{\mu_p}&=\sum_{j=1}^{+\infty}\sum_{p=a_j}^{a_{j+1}-1}\frac{1}{\mu_p}=\sum_{j=1}^{j_0}\sum_{p=a_j}^{a_{j+1}-1}\frac{1}{\mu_p}+\sum_{j=j_0+1}^{+\infty}\sum_{p=a_j}^{a_{j+1}-1}\frac{1}{\mu_p}=\sum_{j=1}^{j_0}\sum_{p=a_j}^{a_{j+1}-1}\frac{1}{\mu_p}+\sum_{j=j_0+1}^{+\infty}\frac{(a_{j+1}-a_j)b_j}{a_j^s}
\\&
\le\sum_{j=1}^{j_0}\sum_{p=a_j}^{a_{j+1}-1}\frac{1}{\mu_p}+\sum_{j=j_0+1}^{+\infty}\frac{1}{2^j}\le \sum_{j=1}^{j_0}\sum_{p=a_j}^{a_{j+1}-1}\frac{1}{\mu_p}+1<+\infty.
\end{align*}

$(c)$ Note that the matrix $\mathcal{M}_{\omega_{\mathbf{M}}}$ is constant if and only if $\mathbf{M}$ satisfies moderate growth; recall Lemma \ref{assoweightomega0} and $(iii)$ in Section \ref{weightfctsection}.

So let us prove that $\mathbf{M}$ violates moderate growth and this is equivalent to having $\sup_{p\in\NN}\frac{\mu_{2p}}{\mu_p}=+\infty$; see again e.g. \cite[Thm. 3.1]{modgrowthstrange} and the references there. By $(iv)$ we get that for all $j\in\NN_{>0}$ sufficiently large it holds that $(a_{j+1}-1<)2(a_{j+1}-1)\le a_{j+2}-1$ and so with the choice $p=a_{j+1}-1$ one has $\frac{\mu_{2(a_{j+1}-1)}}{\mu_{a_{j+1}-1}}=\frac{a_{j+1}^s/b_{j+1}}{a_j^s/b_j}$ and the conclusion follows by $(v)$.\vspace{6pt}

$(d)$ We use some properties mentioned in the proof of Theorem \ref{Thm12igeneral3}. Note that $\mathbf{M},\mathbf{G}^s\in\hyperlink{LCset}{\mathcal{LC}}$ and so both matrices $\mathcal{M}_{\omega_{\mathbf{M}}}$ and $\mathcal{M}_{\omega_{\mathbf{G}^s}}:=\{(\mathbf{G}^s)^{(\iota)}: \iota>0\}$ are standard log-convex; recall Lemma \ref{assoweightomega0} and $(i)$ in Section \ref{weightfctsection}.

First, assume $\mathbf{M}\hyperlink{mtriangle}{\vartriangleleft}\mathbf{G}^s$. Then, in view of Lemma \ref{inclusioncharacterizationlemma} this relation yields that for each $h>0$ there exists $C_h\ge 1$ such that $\omega_{\mathbf{G}^s}(t)\le\omega_{\mathbf{M}}(ht)+C_h$ for all $t\ge 0$ (see \eqref{inclusioncharacterizationlemmaequ1}); this is precisely relation $\omega_{\mathbf{M}}\vartriangleleft_{\mathfrak{c}}\omega_{\mathbf{G}^s}$ in \cite[Sect. 2.3]{ultradifferentiablecomparison}. As commented before, $\omega_{\mathbf{G}^s}$ does have \hyperlink{om6}{$(\omega_6)$} and hence \cite[Prop. 2.2 $(iv)$]{ultradifferentiablecomparison} gives $\omega_{\mathbf{M}}\hyperlink{omvartriangle}{\vartriangleleft}\omega_{\mathbf{G}^s}$; i.e. $t^{1/s}=o(\omega_{\mathbf{M}}(t))$ as $t\rightarrow+\infty$.

Conversely, $\omega_{\mathbf{M}}\hyperlink{omvartriangle}{\vartriangleleft}\omega_{\mathbf{G}^s}$, \cite[Prop. 2.2 $(iii)$]{ultradifferentiablecomparison} and the fact that $\omega_{\mathbf{G}^s}$ satisfies \hyperlink{om6}{$(\omega_1)$} give $\omega_{\mathbf{M}}\vartriangleleft_{\mathfrak{c}}\omega_{\mathbf{G}^s}$. And $\omega_{\mathbf{M}}\vartriangleleft_{\mathfrak{c}}\omega_{\mathbf{G}^s}$ implies $\mathbf{M}\hyperlink{mtriangle}{\vartriangleleft}\mathbf{G}^s$ by $(ii)$ in Lemma \ref{inclusioncharacterizationlemma} since $\mathbf{M}$ is log-convex.

By \cite[Lemma 5.8 $(i)$\&$(iii)$]{genLegendreconjBMT} applied to these associated weight functions, and since $\omega_{\mathbf{G}^s}$ enjoys \hyperlink{om1}{$(\omega_1)$}, we infer that $\omega_{\mathbf{M}}\hyperlink{omvartriangle}{\vartriangleleft}\omega_{\mathbf{G}^s}$ is equivalent to $\mathcal{M}_{\omega_{\mathbf{M}}}\vartriangleleft\mathcal{M}_{\omega_{\mathbf{G}^s}}$. Since $\mathcal{M}_{\omega_{\mathbf{G}^s}}$ is constant and $\mathbf{G}^s=(\mathbf{G}^{s})^{(1)}$, the desired relation \eqref{counterexampleequ} is equivalent to $\mathcal{M}_{\omega_{\mathbf{M}}}\vartriangleleft\mathcal{M}_{\omega_{\mathbf{G}^s}}$. (Indeed, all arguments work when $\mathbf{G}^s$ is replaced by any other sequence $\mathbf{N}\in\hyperlink{LCset}{\mathcal{LC}}$ satisfying moderate growth and having \cite[$(3.2)$]{subaddlike}.)\vspace{6pt}

It remains to show $\mathbf{M}\hyperlink{mtriangle}{\vartriangleleft}\mathbf{G}^s$: The corresponding sequence of quotients for $\mathbf{G}^s$ is given by $(p^s)_{p\in\NN_{>0}}$ and so, by $(ii)$, we get $\frac{a_j^s}{\mu_{a_j}}=b_j\rightarrow+\infty$ as $j\rightarrow+\infty$. Since $\mu_p=\mu_{a_j}$ and $a_j^s\le p^s$ for all $a_j\le p\le a_{j+1}-1$ (see \eqref{exampledef}) it follows that $\lim_{p\rightarrow+\infty}\frac{\mu_p}{p^s}=0$.

Finally, since $\mathbf{M}$ is log-convex with $M_0=1$ and since $\mathbf{G}^s$ has moderate growth, we have $(M_p)^{1/p}\le\mu_p$ and $p^s\le A(p!^s)^{1/p}$ for some $A\ge 1$ and all $p\in\NN_{>0}$ (see again e.g. \cite[Thm. 3.1]{modgrowthstrange}). Summarizing,
$$\exists\;A\ge 1\;\forall\;\epsilon>0\;\exists\;p_{\epsilon}\in\NN_{>0}\;\forall\;p\ge p_{\epsilon}:\;\;\;(M_p)^{1/p}\le\mu_p\le\epsilon p^s\le A\epsilon(p!^s)^{1/p},$$
and so $\mathbf{M}\hyperlink{mtriangle}{\vartriangleleft}\mathbf{G}^s$ holds.\vspace{6pt}

$(e)$ This follows immediately by taking into account $(ii)$ in Section \ref{weightfctsection} for the matrix $\mathcal{M}_{\omega_{\mathbf{M}}}$; recall Lemma \ref{assoweightomega0} and $(a)$.
\qed\enddemo

\bibliographystyle{plain}
\bibliography{Bibliography}
\end{document}